\documentclass[11pt,a4paper,twoside]{article}
\usepackage{graphics}
\usepackage{amssymb}
\usepackage{amsmath}
\usepackage{mathrsfs} %引用花体字母宏包
\usepackage{makeidx}
\usepackage{color}
\usepackage{graphicx}
\usepackage{subfigure}
\usepackage{multicol}
\usepackage{multirow}
\usepackage{float}
\usepackage{booktabs}
\usepackage[colorlinks,
            linkcolor=red,
            anchorcolor=blue,
            citecolor=blue
            ]{hyperref} %引用超链接宏包
\usepackage{caption}
\usepackage{algorithm}
\usepackage{algorithmic}
\usepackage{latexsym, cite, bm, amsthm, amsmath}

\def\appendix{\textbf{Appendix~A.}}

\def\vextra{\vphantom{\vrule height0.4cm width0.9pt depth0.1cm}}

\newtheorem{thm}{Theorem}[section]
\newtheorem{cor}{Corollary}[section]

\newtheorem{prop}{Proposition}[section]
\newtheorem{alg}{Algorithm}[section]
\newtheorem{rem}{Remark}[section]
\newtheorem{exam}{Example}[section]
\numberwithin{equation}{section}

\title{Optimal parameter for the SOR-like iteration method on solving the system
of absolute value equations}

\usepackage[marginal]{footmisc}
\usepackage{authblk}
\author[a]{Cairong Chen\thanks{Supported by the NSFC grants (11901024 and 11672074) and China Postdoctoral Science Foundation (2019M660385). (Email address: cairongchen@buaa.edu.cn).}}
\author[a,b]{Dongmei Yu\thanks{Supported by China Postdoctoral Science Foundation (2019M650449) and Natural Science Foundation of Liaoning Province(Nos. 1584942734251, LJ2019JL005). Corresponding author: D.-M. Yu, (yudongmei1113@163.com).}}
\author[a,c]{Deren Han\thanks{Supported by the NSFC grants 11625105, 11926358 and  11431002. (Email address: handr@buaa.edu.cn).}}
\affil[a]{School of Mathematical Sciences, Beihang University, Beijing, 100191, P.R. China.}
\affil[b]{Institute for Optimization and Decision Analytics, Liaoning Technical University, Fuxin, 123000, P.R. China.}
 \affil[c]{Beijing Advanced Innovation Center for Big Data and Brain Computing (BDBC), Beihang University, Beijing, 100191, P.R. China.}

\begin{document}
\date{\today}
\maketitle

%\double

%\begin{center}
%%%%%%%%%%%%%%%%%%%%%%%%%摘要部分%%%%%%%%%%%%%%%%%%%%%%%%
\begin{quote}
{\bf Abstract:} The SOR-like iteration method for solving
the absolute value equations~(AVE) of finding a vector $x$ such that $Ax - |x| - b = 0$ with $\nu = \|A^{-1}\|_2 < 1$ is investigated. The convergence conditions of the SOR-like iteration method proposed by Ke and Ma ([{\em Appl. Math. Comput.}, 311:195--202, 2017]) are revisited and a new proof is given, which exhibits some  insights in determining the convergent region and the optimal iteration parameter. Along this line, the optimal parameter which minimizes $\|T_\nu(\omega)\|_2$ with
$$T_\nu(\omega) = \left(\begin{array}{cc} |1-\omega| & \omega^2\nu \\  |1-\omega|  &  |1-\omega| +\omega^2\nu \end{array}\right)$$
and the approximate optimal parameter which minimizes $\eta_{\nu}(\omega)
=\max\{|1-\omega|,\nu\omega^2\}$
are explored. The optimal and approximate optimal parameters are iteration-independent and the bigger value of $\nu$ is, the smaller convergent region of the iteration parameter $\omega$ is. Numerical results are presented to demonstrate that the SOR-like iteration method with the optimal parameter is superior to that with the approximate optimal parameter proposed by Guo, Wu and Li ([{\em Appl. Math. Lett.}, 97:107--113, 2019]). In some situation, the SOR-like itration method with the optimal parameter performs better, in terms of CPU time, than the generalized Newton method (Mangasarian, [{\em Optim. Lett.}, 3:101--108, 2009]) for solving the AVE.

{\small

}

\end{quote}

%%%% Start %%%%%%
\section{Introduction}\label{sec:sec1}
\qquad In this paper, we consider the solution of the system of absolute value equations (AVE)
\begin{equation}\label{eq:ave}
    Ax-|x|-b=0,
\end{equation}
where $A\in \mathbb{R}^{n\times n}, b\in \mathbb{R}^n$, and $|x|\in \mathbb{R}^n$ denotes the component-wise absolute value of the vector $x \in \mathbb{R}^n$. The AVE \eqref{eq:ave} can be regarded as a special case of the more general AVE $Ax+B|x|-b=0$ with $B\in \mathbb{R}^{n\times n}$, which was introduced in \cite{rohn2004} and further investigated in  \cite{hlad2018,mang2007,prok2009,wuli2019,wacc2019}; just list a few. The AVE is equivalent to linear complementarity problem and mixed integer programming \cite{huhu2010,mame2006,mang2007,prok2009}, and it recently receives much attention in the optimization community.

It was known that determining the existence of a solution to the AVE is NP-hard \cite{mang2007} and if solvable, the problem of checking whether the AVE has unique or multiple solutions is NP-complete \cite{prok2009}. In addition, one sufficient condition for the AVE \eqref{eq:ave} being uniquely solvable for any $b$ is described in the following proposition.

\begin{prop}[\!\!\cite{mame2006}]\label{prop1} Assume that $A \in \mathbb{R}^{n\times n}$ is nonsingular. If $\|A^{-1}\|_2<1$, then the AVE~\eqref{eq:ave} has a unique solution $x^*$ for any $b \in \mathbb{R}^{n}$.
\end{prop}

Hereafter, we assume that  $A \in \mathbb{R}^{n\times n}$ and $\|A^{-1}\|_2<1$.

Recently, there has been a surge of interest in solving the AVE \eqref{eq:ave} with $\|A^{-1}\|_2<1$ and a large number of numerical iteration methods have been proposed, including the generalized Newton method \cite{mang2009,lilw2018}, the smoothing Newton method \cite{caqz2011}, the Levenberg-Marquardt method~\cite{iqia2015}, the SOR-like iteration method \cite{kema2017,guwl2019}, the generalization of the Gauss-Seidel iteration method \cite{edhs2017} and others; see \cite{ke2020,noin2012,crfp2016,mang2007a,sayc2018,wacl2017,maee2017,abhm2018,maer2018,miys2017,yu2020} and references therein.

Our work here is inspired by recent studies on the SOR-like iteration method for solving the AVE \eqref{eq:ave} \cite{kema2017,guwl2019}. The SOR-like iteration method is one-parameter-dependent and thus it is an important problem to determine the optimal iteration parameter. By an optimal iteration parameter, we mean it is the iteration parameter such that the SOR-like iteration method gets the fastest convergence rate. Usually, it seems not to be an easy task to find the optimal value of the involved iteration parameter; while it remains significance to find a somewhat optimal one. Recently, Guo, Wu and Li \cite{guwl2019} obtained an optimal iteration parameter that minimizes the spectral radius of the iteration matrix. However, it is associated with the spectral radius $\rho(D(x^{(k+1)})A^{-1})$ and thus it is iteration-dependent (that is, it varies with the iteration sequence $\{x^{(k)}\}$). To compute the optimal parameter in every iterative step is expensive, especially when the dimension of $A$ is large, thus an approximate one is used in the numerical experiments of \cite{guwl2019}.

The goal of this paper is twofold: to revisit the convergence conditions of the SOR-like iteration method on solving the AVE \eqref{eq:ave} and to explore respectively an optimal iteration parameter and an approximate optimal iteration parameter for the the SOR-like iteration method. It is important that our optimal and approximate optimal parameters are iteration-independent. Numerical results demonstrate that the SOR-like iteration method with our optimal iteration parameter is superior to that with the approximate one proposed in \cite{guwl2019} for solving the AVE~\eqref{eq:ave}. The SOR-like iteration method with our approximate optimal iteration parameter is better than that with the approximate one proposed in \cite{guwl2019} for solving the AVE~\eqref{eq:ave} in some cases. In addition, the SOR-like itration method with the optimal parameter sometimes performs better, in terms of CPU time, than the generalized Newton method~\cite{mang2009}.

The rest of this paper is organized as follows. In Section~\ref{sec:cond} we revisit the convergence conditions of the SOR-like iteration method for solving the AVE \eqref{eq:ave}. Section~\ref{sec:opt} characterizes the optimal and approximate optimal iteration parameters for the SOR-like iteration method. In Section~\ref{sec:numer}, some numerical examples are given to demonstrate our claims made in the previous sections. Finally, some concluding remarks are given in Section~\ref{sec:con}.

{\bf Notations.}
$\mathbb{R}^{n\times n}$ is the set of all $n \times n$ real matrices and $\mathbb{R}^{n}= \mathbb{R}^{n\times 1}$. $I$ is the identity matrix with suitable dimension. For $x \in \mathbb{R}^{n}$, $x_i$ refers to its $i$th entry, $|x|$ is in $\mathbb{R}^{n}$ with its $i$th entry $| x_i |$, and  $| \cdot |$ denotes absolute value. $\text{sgn}(x)$ denotes a vector with components equal to $-1, 0$ or $1$ depending on whether the corresponding component of the vector $x$ is negative, zero or positive, respectively.
For $x \in \mathbb{R}^{n}$, $\text{diag}(x)$ represents a diagonal matrix with $x_i$ as its diagonal entries for every $i=1,2,\cdots,n$. $\text{tridiag}(a, b, c)$ denotes a tridiagonal matrix that has $a, b, c$ as the subdiagonal, main diagonal and superdiagonal entries in the matrix, respectively; $\text{Tridiag}(A, B, C)$ denotes a block tridiagonal matrix that has $A, B, C$ as the subdiagonal, main diagonal and superdiagonal block entries in the matrix, respectively. $\|A\|_2$ denotes the spectral norm of $A$ and is defined by the formula $\| A \|_2\doteq \max \left\{ \| A x \|_2 : x \in \mathbb{R}^{n}, \|x\|_2=1 \right\}$, where $\| x \|_2 $ is the $2$-norm of vector $x \in \mathbb{R}^{n}$. When $X$ is a square matrix, we denote by $\rho(X)$, $\lambda_{\max}(X)$ and $\lambda_{\min}(X)$ its spectral radius, largest eigenvalue and smallest eigenvalue, respectively. For a function $f$, $f'$ denotes its first derivative. $\lim_{x\rightarrow x_0^+}f(x)$ and $\lim_{x\rightarrow x_0^-}f(x)$ denote the right and left limits of a function $f(x)$ at a point~$x_0$, respectively. $\textbf{1}_n$ denotes the $n$-vector of all ones.

\section{Revisit the convergence conditions of the SOR-like iteration method}\label{sec:cond}
\qquad The SOR-like iteration method for solving the AVE \eqref{eq:ave} is firstly proposed by Ke and Ma \cite{kema2017}. After that, in a different perspective, Guo, Wu and Li in \cite{guwl2019} present some new convergence conditions of the SOR-like iteration method. In this section,  along the lines in \cite{kema2017}, we will further study the convergence conditions of the SOR-like iteration method for solving the AVE~\eqref{eq:ave}. For this purpose, we first briefly review the SOR-like iteration method for solving the AVE~\eqref{eq:ave}.

The AVE \eqref{eq:ave} is equivalent to
\begin{equation}\label{eq:equiv}
\left\{\begin{array}{l}
Ax-y=b,\\
y-|x|=0,
\end{array}\right.
\end{equation}
that is,
\begin{align*}
\boldsymbol{Az}\doteq
\left(\begin{array}{cc}
A    &   -I \\
-D(x) &   I
\end{array}\right)
\left(\begin{array}{c}
x    \\   y
\end{array}\right)
=
\left(\begin{array}{c}
b    \\  0
\end{array}\right)
\doteq\boldsymbol b,
\end{align*}
where $D(x)\doteq\text{diag}(\text{sgn}(x))$. By splitting the coefficient matrix $\boldsymbol{A}= \mathcal{D} - \mathcal{L} - \mathcal{U} $ with
\begin{align*}
\mathcal{D}=
\left(\begin{array}{cc}
 A   &   0 \\
 0   &   I
\end{array}\right),\quad
\mathcal{L}=
\left(\begin{array}{cc}
0     &   0 \\
D(x)   &   0
\end{array}\right),\quad
\mathcal{U}=
\left(\begin{array}{cc}
0      &   I \\
0      &   0
\end{array}\right),
\end{align*}
one can obtain that the iteration scheme of the SOR-like iteration method is
\begin{align*}
\left(\begin{array}{c}
 x^{(k+1)}   \\   y^{(k+1)}
\end{array}\right)
=\mathcal{M}_\omega
\left(\begin{array}{c}
 x^{(k)}   \\   y^{(k)}
\end{array}\right)
+\omega(\mathcal{D} - \omega \mathcal{L})^{-1}
\left(\begin{array}{c}
 b    \\  0
\end{array}\right),
\end{align*}
where $\mathcal{M}_\omega = (\mathcal{D} - \omega \mathcal{L})^{-1}\left[(1-\omega) \mathcal{D}+\omega\mathcal{U}\right]$, and $\omega>0$ is the iteration parameter. Specifically, the SOR-like iteration method is described in the following Algorithm~\ref{Alg1}.

\begin{alg}[\!\!\cite{kema2017}]\label{Alg1}$($SOR-like iteration method for solving the AVE \eqref{eq:ave}$)$
Let $A\in \mathbb{R}^{n\times n}$ be a nonsingular matrix and $b\in \mathbb{R}^{n}$. Given initial vectors $x^{(0)}\in \mathbb{R}^{n}$ and $y^{(0)}\in \mathbb{R}^{n}$, for $k=0,1,2,\cdots$ until the iteration sequence $\left\{(x^{(k)},y^{(k)})\right\}_{k=0}^\infty$ is convergent, compute
\begin{align}\label{eq:sorl2}
\left\{\begin{array}{l}
x^{(k+1)}=(1-\omega)x^{(k)}+\omega A^{-1}(y^{(k)}+b),\\
y^{(k+1)}=(1-\omega)y^{(k)}+\omega |x^{(k+1)}|,
\end{array}\right.
\end{align}
where the iteration parameter $\omega $ is a positive constant.
\end{alg}

Now we are in position to explore the convergence of the SOR-like iteration method for solving the AVE \eqref{eq:ave}. Let $(x^*,y^*)$ be the solution pair of the nonlinear equation \eqref{eq:equiv} and denote $e_k^x=x^*-x^{(k)}$ and $e_k^y=y^*-y^{(k)}$, where $(x^{(k)},y^{(k)})$ is generated by Algorithm~\ref{Alg1}. Then, for the SOR-like iteration method, the following results are known.

\begin{thm}[\!\!\cite{kema2017}]\label{thm1} Assume that $A \in \mathbb{R}^{n\times n}$ is a nonsingular matrix and $b\in \mathbb{R}^{n}$. Denote
$$ \nu=\|A^{-1}\|_2, \quad a=|1-\omega|\quad \text{and}\quad c=\omega^2\nu. $$
If
\begin{equation}\label{eq:cond1}
0<\omega< 2  \qquad \text{and} \qquad  a^4-3a^2 -2ac- 2c^2 +1 >0,
\end{equation}
then the following inequality
\begin{equation}\label{eq:res}
\| |(e_{k+1}^x,e_{k+1}^y)| \|_{\omega} < \| |(e_k^x,e_k^y) |\|_{\omega}
\end{equation}
holds for $ k=0,1,2,\cdots $. Here the norm $\| |\cdot|\|_{\omega}$ is defined by
$$\| |(e_k^x,e_k^y) |\|_{\omega}:=\sqrt {\|e_k^x \|_2^2+\omega ^{-2}\|e_k^y \|_2^2 }.$$
\end{thm}

\begin{cor}[\!\!\cite{kema2017}]\label{cor1}
Let $A \in \mathbb{R}^{n\times n}$ be a nonsingular matrix, $b\in \mathbb{R}^{n}$ and $\nu=\|A^{-1}\|_2$. If
\begin{equation}\label{eq:cond2}
\nu<1\quad\text{and}\quad 1-\tau<\omega<\min\left\{1+\tau,\sqrt{\frac{\tau}{\nu}}\right\},
\end{equation}
where $\tau=\frac{2}{3+\sqrt{5}}$. Then the following inequality
\begin{equation*}
\| |(e_{k+1}^x,e_{k+1}^y)| \|_{\omega} < \| |(e_k^x,e_k^y) |\|_{\omega}
\end{equation*}
holds for $ k=0,1,2,\cdots $.
\end{cor}

The proofs of Theorem~\ref{thm1} and Corollary~\ref{cor1} are given in \cite{kema2017}. However, we will present the following proof of Theorem \ref{thm1} because the expression of the eigenvalue \eqref{eq:lambda} is needed later for discussing the optimal iteration parameter. The former part of the proof is known in \cite{kema2017} and the later part seems new.

{\bf The proof of the Theorem \ref{thm1}:}

\begin{proof} It follows from the proof of Theorem~3.1 in \cite{kema2017} that
\begin{align*}
\left(\begin{array}{c}  \|e_{k+1}^x \|_2  \\   \omega^{-1}\| e_{k+1}^y \|_2   \end{array}\right)
\leq
\left(\begin{array}{cc} |1-\omega| & \omega^2\nu \\  |1-\omega|  &  |1-\omega| +\omega^2\nu \end{array}\right)
\left(\begin{array}{c}  \|e_{k}^x \|_2  \\   \omega^{-1} \| e_{k}^y \|_2   \end{array}\right),
\end{align*}
that is,
\begin{equation*}%\label{eq:iter}
\| |(e_{k+1}^x,e_{k+1}^y)| \|_{\omega} \le \|T_\nu(\omega) \|_2 \cdot \| |(e_k^x,e_k^y) |\|_{\omega}
\end{equation*}
with
$$ T_\nu(\omega)=\left(\begin{array}{cc} a & c \\  a  &  a+c \end{array}\right). $$
In order to prove the inequality \eqref{eq:res}, we turn to consider the choice of the parameter $\omega$ such that $\| T_\nu(\omega) \|_2 <1$.

Since
$$ H_\nu(\omega)=T_\nu( \omega)^T T_\nu( \omega)= \left(\begin{array}{cc}  2a^2 & a^2+2ac  \\ a^2+2ac   & a^2+2c^2+2ac  \end{array}\right)$$
is a symmetric positive semi-definite matrix, we have $\| T_\nu(\omega) \|_2^2 = \rho\left( T_\nu(\omega)^T T_\nu(\omega) \right) =
 \lambda_{\max} \left(H_\nu(\omega)\right)$. Assume that $\lambda$ is an eigenvalue of $H_\nu(\omega)$. Then
$$(\lambda-2a^2)[\lambda-(a^2+2c^2+2a c)]-(a^2+2a c)^2=0,$$
namely,
\begin{equation*}
\lambda^2-(3a^2+2c^2+2a c)\lambda+a^4=0,
\end{equation*}
from which we obtain
$$ \lambda= \frac{3a^2+2c^2+2a c \pm \sqrt{(3a^2+2c^2+2a c)^2-4a^4}}{2}.$$
Consequently,
\begin{equation}\label{eq:lambda}
\lambda_{\max}(H_\nu(\omega))=\frac{3a^2+2c^2+2a c + \sqrt{(3a^2+2c^2+2a c)^2-4a^4}}{2}.
\end{equation}
In particular,
\begin{align*}
\lambda_{\max}(H_\nu(\omega))<1 &\Longleftrightarrow 3a^2+2c^2+2ac + \sqrt{(3a^2+2c^2+2a c)^2-4a^4}<2\\
&\Longleftrightarrow \sqrt{(3a^2+2c^2+2a c)^2-4a^4}<2-(3a^2+2c^2+2a c).
\end{align*}
Hence, a sufficient condition for the convergence is
\begin{equation}\label{ieq:sc1}
\left\{
\begin{array}{l}
3a^2+2c^2+2ac<2, \\
3a^2+2c^2+2a c<1+a^4.
\end{array}
\right.
\end{equation}
From \eqref{ieq:sc1}, we have $\lambda_{\max}(H_\nu(\omega))<1$ provided that $1+a^4<2$ and $3a^2+2c^2+2a c<1+a^4$. It is easy to check that $1+a^4<2$ is equivalent to $0<\omega<2$, which completes the proof.
\end{proof}

The conditions in~\eqref{eq:cond1} guarantee the convergence of the SOR-like iteration method for solving the AVE~\eqref{eq:ave}. However, when comparing them with the conditions in~\eqref{eq:cond2}, the second inequality in~\eqref{eq:cond1} seems harder to be checked at first glance. Thus, we will go ahead and talk something more about it.

\begin{figure}[t]
{\centering
\begin{tabular}{ccc}
\hspace{-0.3 cm}
\resizebox*{0.50\textwidth}{0.26\textheight}{\includegraphics{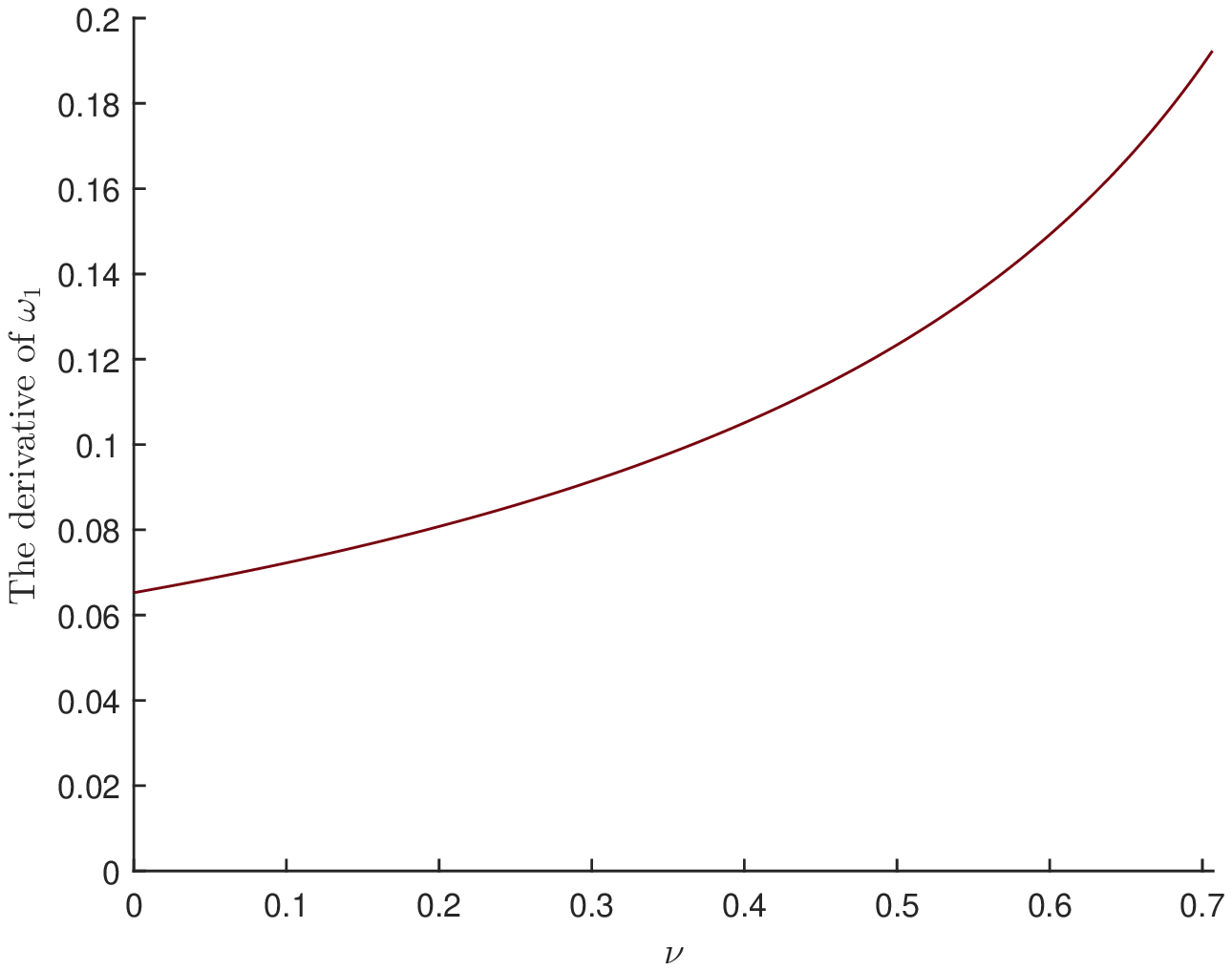}}
& & \hspace{-0.9 cm}
\resizebox*{0.50\textwidth}{0.26\textheight}{\includegraphics{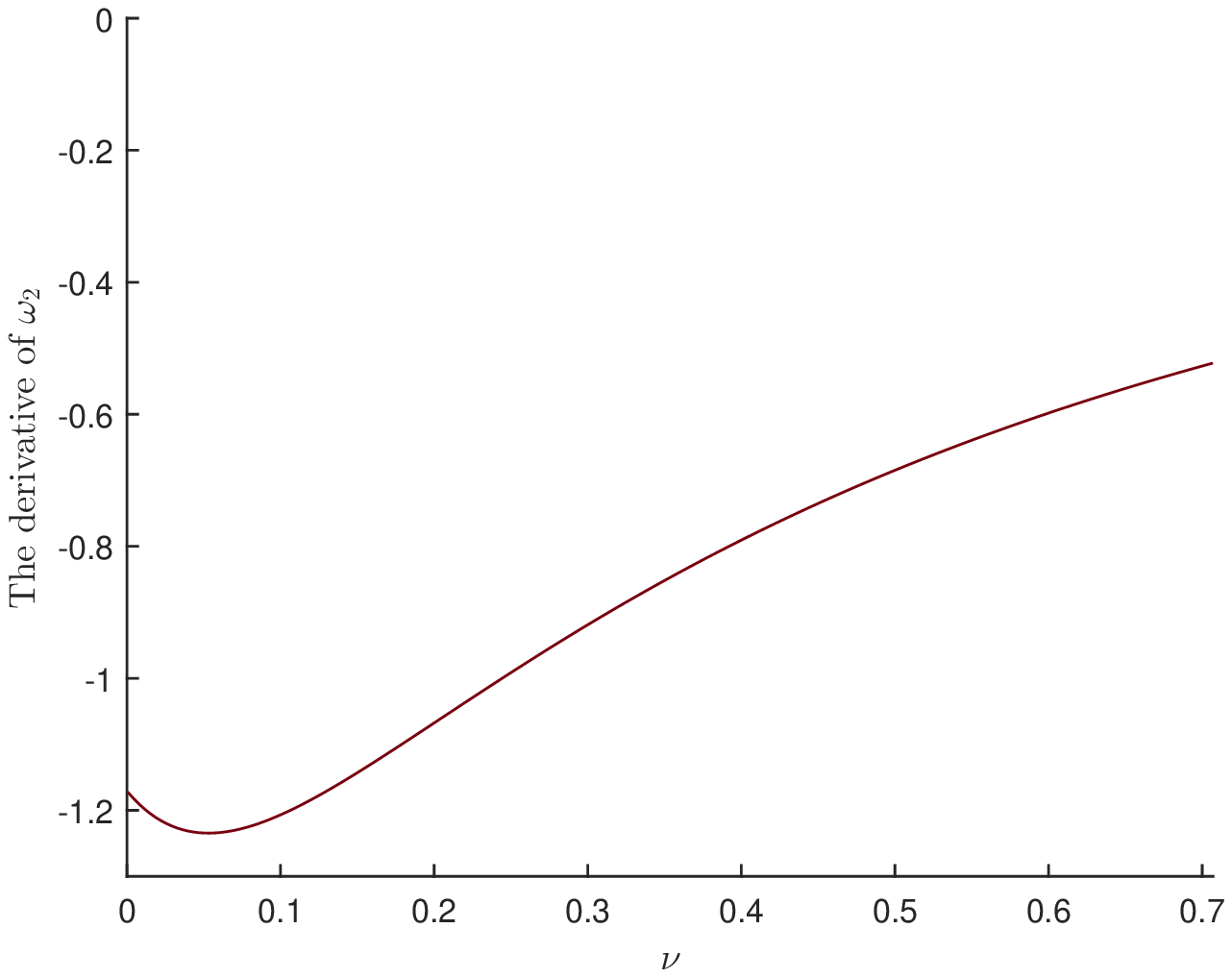}}\vspace{2ex}\\
\hspace{-0.3 cm}
\resizebox*{0.50\textwidth}{0.26\textheight}{\includegraphics{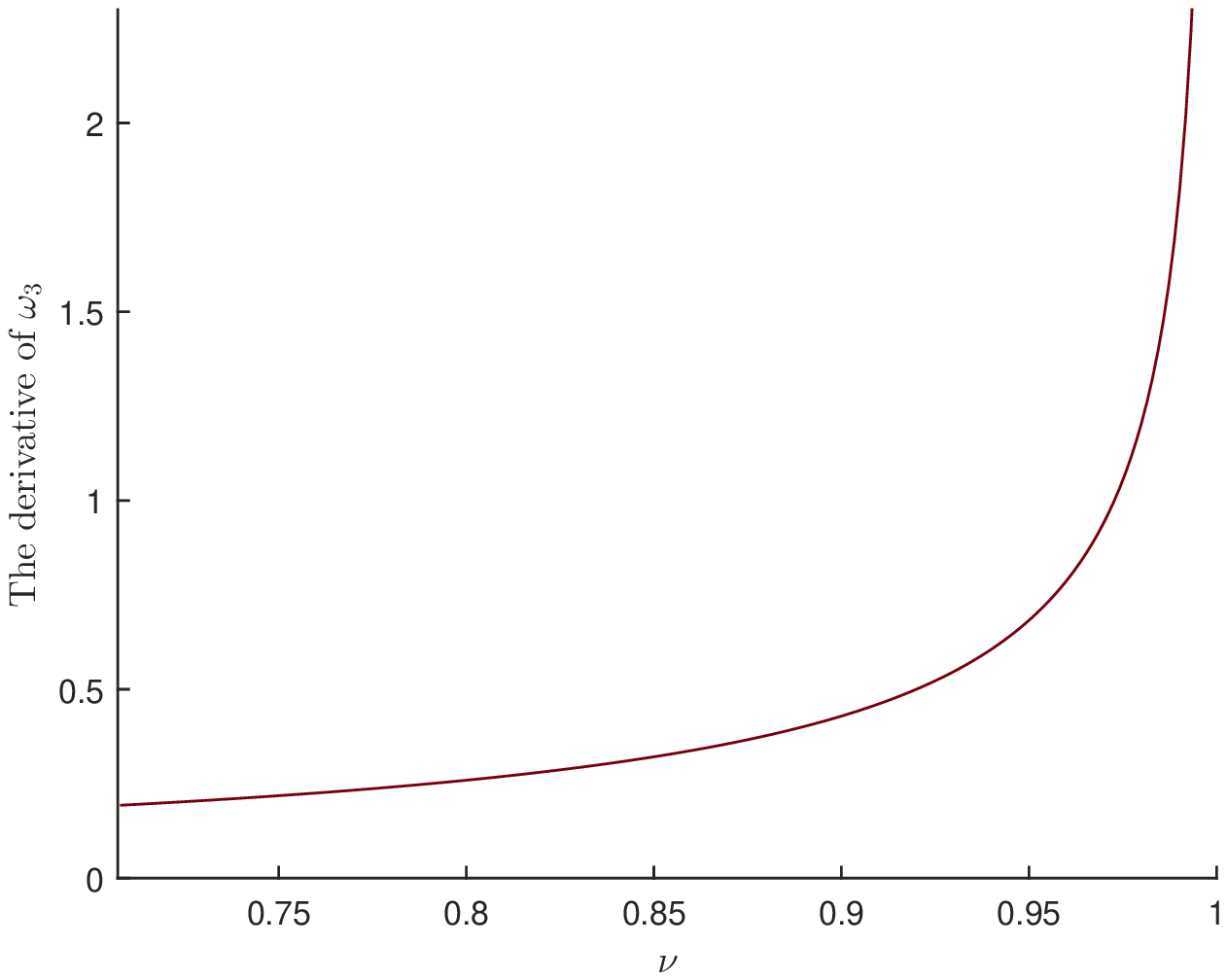}}
& & \hspace{-0.9 cm}
\resizebox*{0.50\textwidth}{0.26\textheight}{\includegraphics{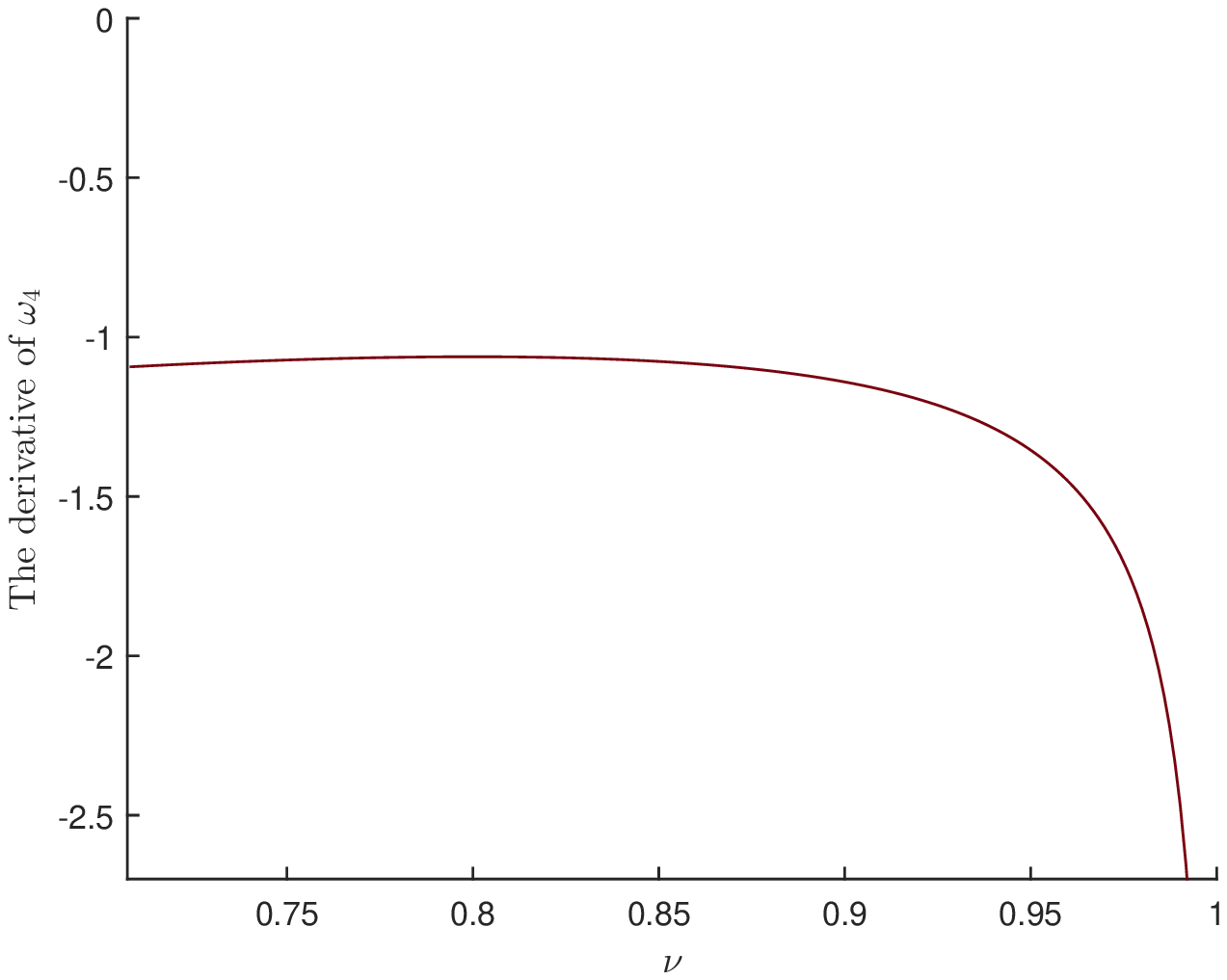}}
\end{tabular}\par
}\vspace{-0.15 cm}
\caption{Curves for derivatives of $\omega_{i}(\nu)$ with $i=1,\,2,\,3,\,4$.}
\label{Fig:domega}
\end{figure}

\begin{figure}[t]
{\centering
\begin{tabular}{ccc}
\hspace{-0.3 cm}
\resizebox*{0.50\textwidth}{0.26\textheight}{\includegraphics{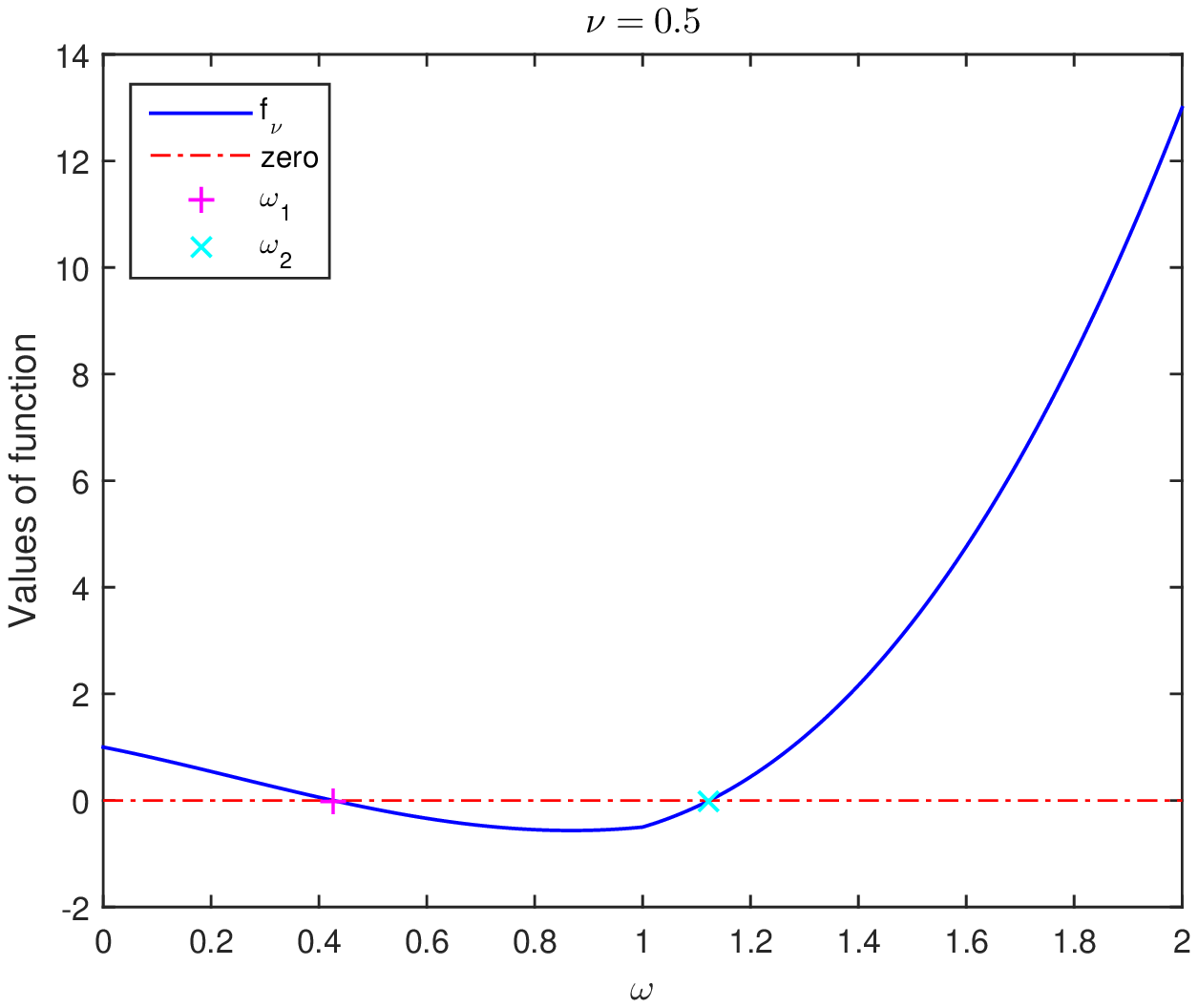}}
& & \hspace{-0.9 cm}
\resizebox*{0.50\textwidth}{0.26\textheight}{\includegraphics{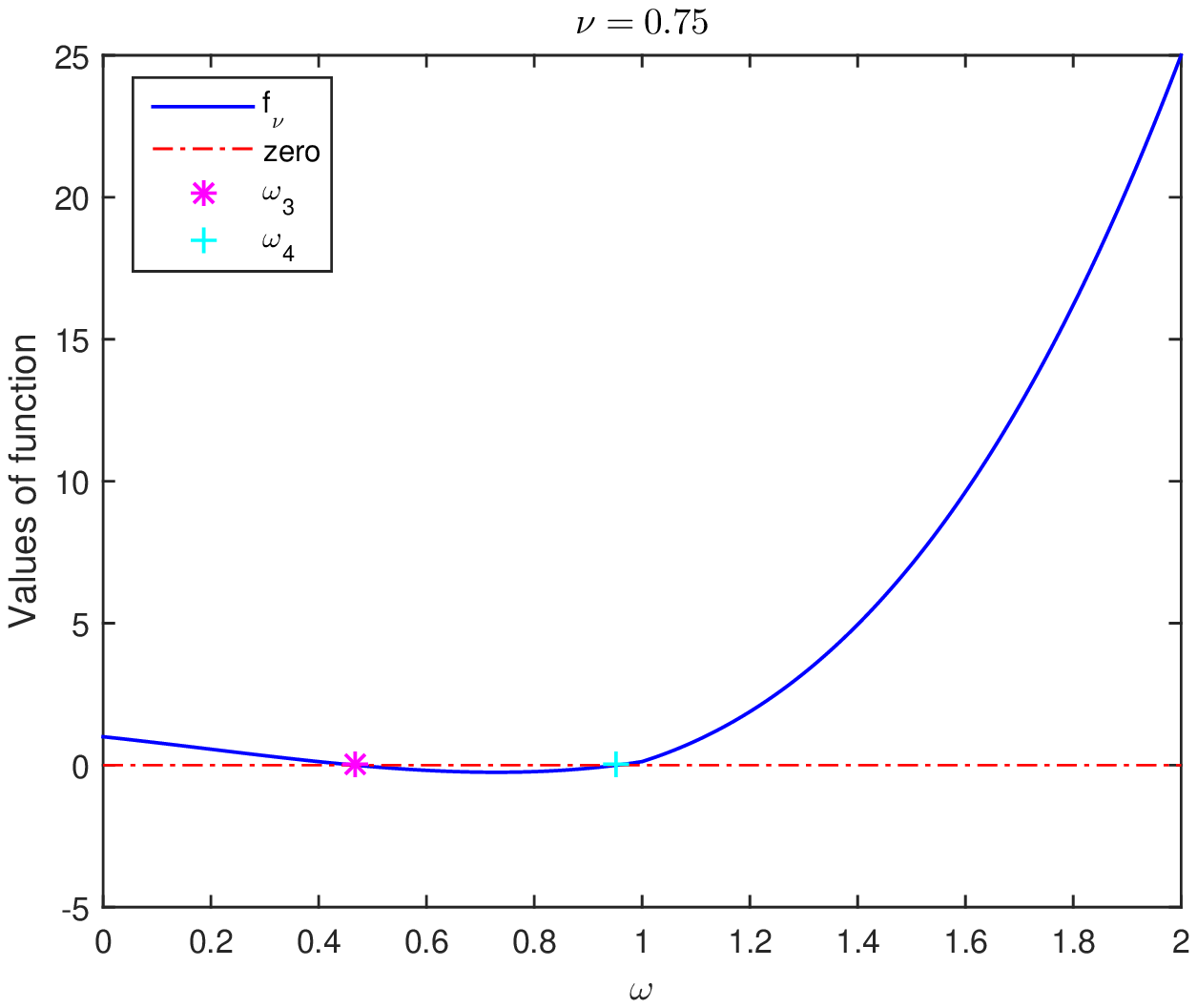}}\vspace{2ex}\\
\hspace{-0.3 cm}
\resizebox*{0.50\textwidth}{0.26\textheight}{\includegraphics{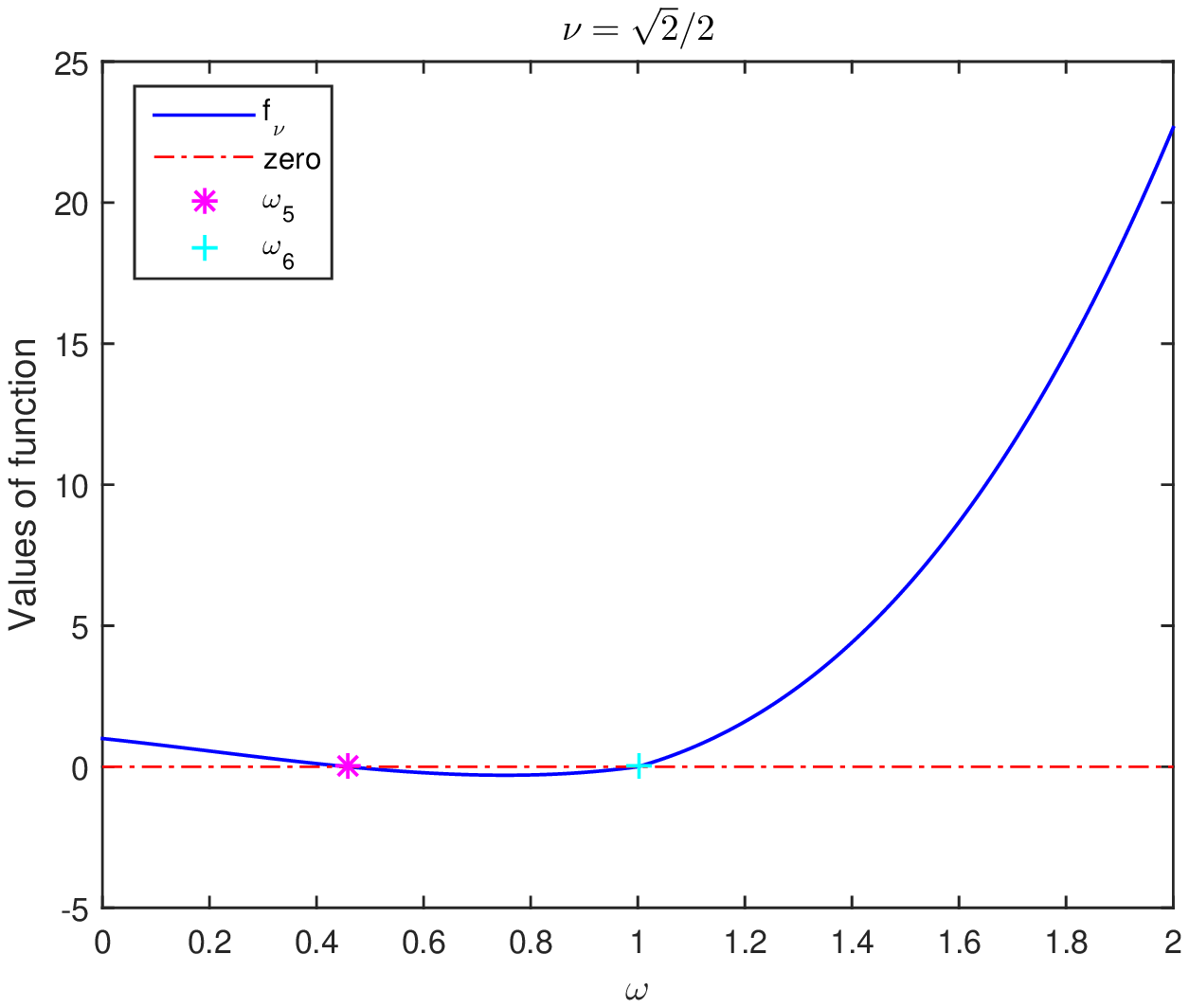}}
& & \hspace{-0.9 cm}
\resizebox*{0.50\textwidth}{0.26\textheight}{\includegraphics{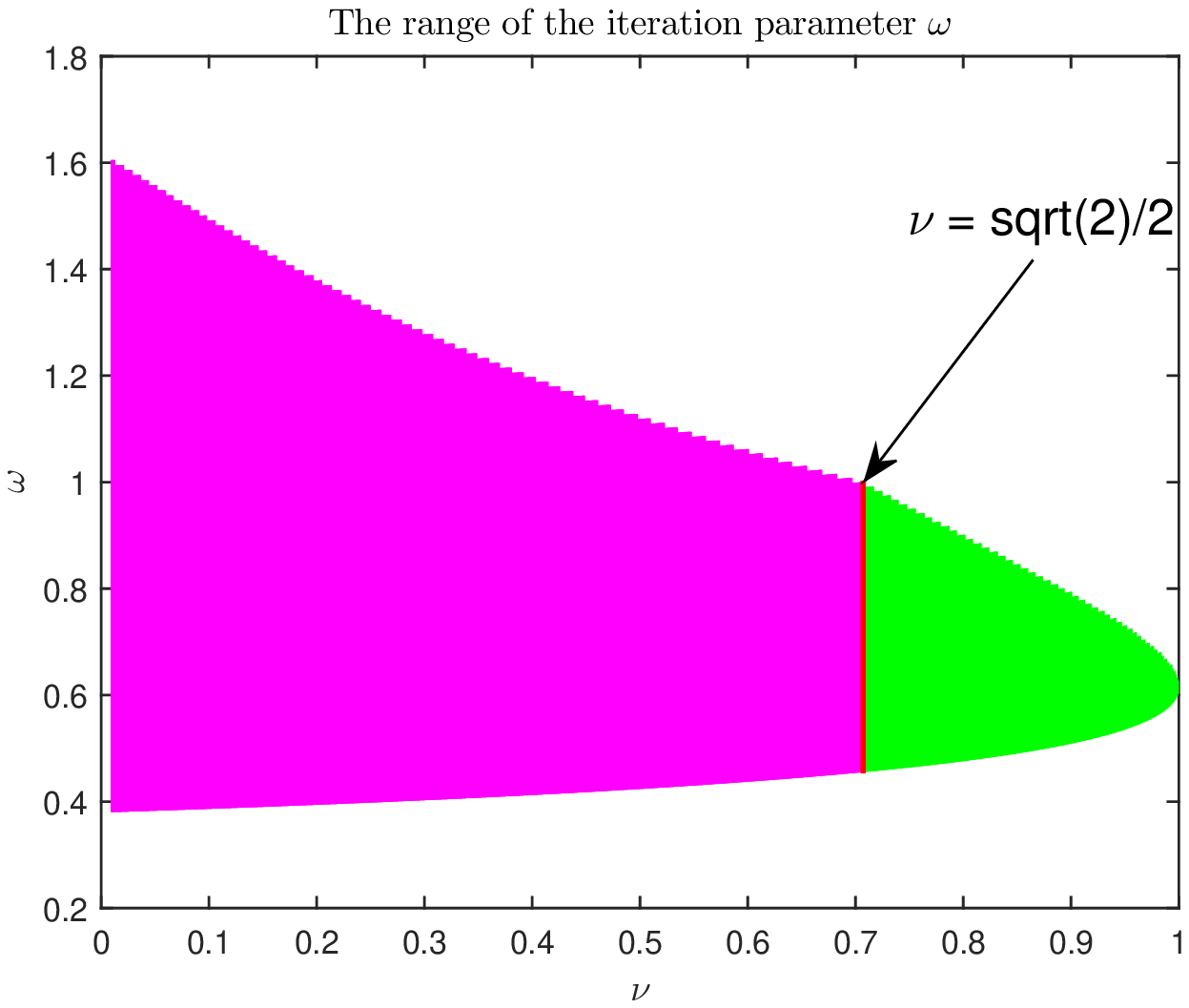}}
\end{tabular}\par
}\vspace{-0.15 cm}
\caption{Curves for function $f_{\nu}$ with $\nu=0.5$ (top-left), $\nu=0.75$ (top-right), $\nu=\frac{\sqrt{2}}{2}$ (bottom-left) and the range of the iteration parameter $\omega$ (bottom-right).}
\label{Fig:vof}
\end{figure}

Let
\begin{align}\nonumber
f_{\nu}(\omega)& = 3a^2+2c^2+2a c - a^4 -1\\\label{eq:f}
         & = 3(\omega - 1)^2 + 2\nu^2\omega^4 + 2\nu\omega^2|\omega - 1| - (\omega - 1)^4 - 1.
\end{align}
According to Proposition~\ref{prop1} and Theorem~\ref{thm1}, the problem to find $\omega\in (0,2)$ such that $f_{\nu}(\omega)<0$ with $0<\nu<1$ is of interest. For this purpose, by using the computerized algebra system Maple, we first compute the real zeros of the function $f_{\nu}$ defined as in \eqref{eq:f} with respect to $\omega$ ($\nu$ is seen as a parameter). The roots of $f_{\nu}(\omega) = 0$ which locate in $(0,2)$ for $\nu \in (0,1)$ are (some functions of $\nu$):
\begin{align*}
\text{If}\quad 0<\nu<\frac{\sqrt{2}}{2}: &\quad\begin{array}{l}
  \omega_1(\nu) = -\frac{\alpha}{4\beta}+\frac{\gamma}{2\beta}-\frac{\sqrt{-(8\nu^3 - 16\nu^2 + 4\nu -1) -( 8\nu^2 - 2\nu )\gamma}}{2\beta} \\
  \omega_2(\nu) = -\frac{\xi}{4\beta}-\frac{\zeta}{2\beta}+\frac{\sqrt{(8\nu^3 + 16\nu^2 + 4\nu +1) + (8\nu^2 + 2\nu)\zeta}}{2\beta}
\end{array}, \\
\text{If}\quad \frac{\sqrt{2}}{2}<\nu<1: &\quad\begin{array}{l}
  \omega_3(\nu) = -\frac{\alpha}{4\beta}+\frac{\gamma}{2\beta}+\frac{\sqrt{-(8\nu^3 - 16\nu^2 + 4\nu -1) - (8\nu^2 - 2\nu)\gamma}}{2\beta}\\
  \omega_4(\nu) = -\frac{\alpha}{4\beta}-\frac{\gamma}{2\beta}+\frac{\sqrt{-(8\nu^3 - 16\nu^2 + 4\nu -1)+(8\nu^2 - 2\nu)\gamma}}{2\beta}
\end{array},\\
\text{If}\quad \nu = \frac{\sqrt{2}}{2}: &\quad\begin{array}{l}
  \omega_5 =-\frac{1}{7}-\frac{\sqrt {2}}{28}+\frac{\sqrt {242+64\sqrt {2}}}{28} \approx 0.4579\\
  \omega_6 =1
\end{array},
\end{align*}
where
\begin{align*}
\alpha &= 4-2\nu >0,\quad \beta=2\nu^2-1,\quad\xi = 4+2\nu >0,\\
\gamma&=\sqrt{-(\nu-1)(\nu+5)}>0,\quad \zeta = \sqrt{-(\nu+1)(\nu-5)}>0.
\end{align*}
Apparently, we have $0 < \omega_5 <\omega_6 <2$. In Appendex~A, we prove that $\omega_i(\nu)\,(i=1,\,2,\,3,\,4)$ are real functions (with respect to $\nu$) in their domains, respectively. In addition, the curves of the derivatives of $\omega_1$, $\omega_2$ in $\nu \in (0, \frac{\sqrt{2}}{2})$ and $\omega_3$, $\omega_4$ in $\nu \in (\frac{\sqrt{2}}{2}, 1)$ are shown in Figure~\ref{Fig:domega}, from which we can find that $\omega_1$ and $\omega_3$ are strictly monotonously increasing in $\nu\in(0,\frac{\sqrt{2}}{2})$ and $\nu \in (\frac{\sqrt{2}}{2}, 1)$, resepctively, while $\omega_2$  and $\omega_4$ are strictly monotonously decreasing in $\nu \in (0, \frac{\sqrt{2}}{2})$ and $\nu \in (\frac{\sqrt{2}}{2}, 1)$, respectively. Furthermore, by using the L'Hospital's rule once or twice if needed, we have
\begin{align}\label{eq:limits1}
&\lim_{\nu\rightarrow 0^+}\omega_1(\nu) = \frac{3-\sqrt{5}}{2},\; \lim_{\nu\rightarrow \frac{\sqrt{2}}{2}^-}\omega_1(\nu) = \omega_5,\;
\lim_{\nu\rightarrow 0^+}\omega_2(\nu) = \frac{\sqrt{5}+1}{2},\; \lim_{\nu\rightarrow \frac{\sqrt{2}}{2}^-}\omega_2(\nu) = 1,\\\label{eq:limits2}
&\lim_{\nu\rightarrow \frac{\sqrt{2}}{2}^+}\omega_3(\nu) = \omega_5,\; \lim_{\nu\rightarrow 1^-}\omega_3(\nu) = \frac{\sqrt{5}-1}{2}, \; \lim_{\nu\rightarrow \frac{\sqrt{2}}{2}^+}\omega_4(\nu) =1,\;
\lim_{\nu\rightarrow 1^-}\omega_4(\nu) = \frac{\sqrt{5}-1}{2}.
\end{align}
It follows from~\eqref{eq:limits1}, \eqref{eq:limits2} and the strictly monotonous properties of $\omega_i(\nu)\,(i=1,\,2,\,3,\,4)$ that $\omega_1(\nu)\in (\frac{3-\sqrt{5}}{2},\omega_5)\subset (0,2)$ and $\omega_2(v)\in (1,\frac{\sqrt{5}+1}{2})\subset (0,2)$ and thus $0<\omega_1(v)<\omega_2(v)<2$ for $\nu\in(0,\frac{\sqrt{2}}{2})$. Similarly, $\omega_3(\nu)\in (\omega_5,\frac{\sqrt{5}-1}{2})\subset (0,2)$ and $\omega_4(v)\in (\frac{\sqrt{5}-1}{2},1)\subset (0,2)$ and thus $0<\omega_3(v)<\omega_4(v)<2$ for $\nu\in(\frac{\sqrt{2}}{2},1)$. Moreover, $f_{\nu}(\omega) < 0$ if one of the following conditions holds:
\begin{align}\label{eq:ncond1}
&\omega\in (\omega_1(\nu),\omega_2(\nu))\doteq\Omega_1,\quad \text{when}\quad\nu\in (0,\frac{\sqrt{2}}{2}),\\\label{eq:ncond2}
&\omega\in (\omega_3(\nu),\omega_4(\nu))\doteq\Omega_2,\quad \text{when}\quad \nu\in (\frac{\sqrt{2}}{2},1),\\\label{eq:ncond3}
&\omega\in (\omega_5,\omega_6)\doteq\Omega_3 ,\quad  \text{when}\quad  \nu=\frac{\sqrt{2}}{2}.
\end{align}
Then it again can be concluded from~\eqref{eq:limits1}, \eqref{eq:limits2} and the the strictly monotonous properties of $\omega_i(\nu)\,(i=1,\,2,\,3,\,4)$ that the bigger value of $\nu$ is, the smaller the range of $\omega$ will be. In addition, since $\lim_{\nu\rightarrow 1^-}\omega_3(\nu) = \lim_{\nu\rightarrow 1^-}\omega_4(\nu) = \frac{\sqrt{5}-1}{2}$, the value of $\omega$ satisfied conditions~\eqref{eq:cond1} must close to $\frac{\sqrt{5}-1}{2}$ as $\nu \rightarrow 1^-$. Figure~\ref{Fig:vof} plots the curves of function $f_{\nu}$ for some values of $\nu$ and the range of the iteration parameter $\omega$, from which our claims are intuitively shown.

Finally, according to the discussion above, we will rewrite Theorem~\ref{thm1} as follows, i.e., Theorem~\ref{nthm1}, which exhibits more information than the original one.
\begin{thm}\label{nthm1} Assume that $A \in \mathbb{R}^{n\times n}$ is a nonsingular matrix and $b\in \mathbb{R}^{n}$. If one of the conditions \eqref{eq:ncond1}-\eqref{eq:ncond3} holds, then the following inequality
\begin{equation}\label{ie:cov}
\| |(e_{k+1}^x,e_{k+1}^y)| \|_{\omega} \le \|T_\nu(\omega)\|_2 \cdot \| |(e_k^x,e_k^y) |\|_{\omega}\quad \text{with}\quad \|T_\nu(\omega)\|_2<1
\end{equation}
holds for $ k=0,1,2,\cdots $, and \eqref{ie:cov} implies that the sequence $\{x^{(k)}\}$ generated by the SOR-like iteration scheme~\eqref{eq:sorl2} converges to the unique solution of the AVE~\eqref{eq:ave}. Here $\| |(e_k^x,e_k^y) |\|_{\omega}$ is defined as in Theorem~\ref{thm1}.
\end{thm}

%%%%%%%%%%%%%%%%%%%%%%%%%%
\section{Optimal iteration parameter for the SOR-like iteration method}
\label{sec:opt}
\qquad As is known, the SOR-like iteration method involves one parameter and to determine the somewhat optimal iteration parameter is an important problem. From \eqref{ie:cov}, one can obtain that
\begin{equation}\label{eq:err}
0\le \| |(e_{k+1}^x,e_{k+1}^y)| \|_{\omega} \le \|T_\nu(\omega) \|_2 \cdot \| |(e_k^x,e_k^y) |\|_{\omega} \le\cdots\le \|T_\nu(\omega) \|_2^{k+1} \cdot \| |(e_0^x,e_0^y) |\|_{\omega}.
\end{equation}
If the conditions of Theorem~\ref{nthm1} hold, then we have $\|T_\nu(\omega) \|_2<1$ and $\lim_{k\rightarrow\infty}\| |(e_k^x,e_k^y) |\|_{\omega}=0$, that is, the sequence $\{x^{(k)}\}$ generated by the SOR-like iteration scheme \eqref{eq:sorl2} converges to the unique solution of the AVE~\eqref{eq:ave}. In addition, from \eqref{eq:err}, the smaller value of $\|T_\nu(\omega) \|_2$ is, the faster the SOR-like iteration method will converge later on. The question is for what $\omega\in (0,2)$, $\|T_\nu(\omega) \|_2$ is minimized with some given $\nu \in (0,1)$. According to the proof of Theorem~\ref{thm1}, it is equivalent to determining the optimal parameter $\omega\in (0,2)$ that minimizes the eigenvalue $\lambda_{\max}(H_\nu(\omega))$ defined as in \eqref{eq:lambda} for a given $\nu \in (0,1)$.

From \eqref{eq:lambda}, let
\begin{align}\nonumber
g_{\nu}(\omega)&= 3a^2+2c^2+2a c + \sqrt{(3a^2+2c^2+2a c)^2-4a^4}\\\label{eq:fg}
&=3(\omega-1)^2 + 2\nu^2\omega^4+2\nu\omega^2|\omega-1|+\sqrt{\left[3(\omega-1)^2 + 2\nu^2\omega^4+2\nu\omega^2|\omega-1|\right]^2-4(\omega-1)^4},
\end{align}
then minimization of $\lambda_{\max}(H_\nu(\omega))$ is equivalent to minimize $g_{\nu}(\omega)$. Notice that, the function $g_{\nu}$ in \eqref{eq:fg} is continuous but non-smooth with $0<\nu<1$ due to the non-differentiability of the absolute value function. Indeed, it is just not differentiable at $\omega=1$. In addition, by simple calculation, we have
\begin{equation}\label{eq:dgg}
 g'_{\nu}(\omega)=
\left\{\begin{array}{l}
s_\nu(\omega)+\frac{r_\nu(\omega)s_\nu(\omega)
-8(\omega-1)^3}{\sqrt{[r_\nu(\omega)]^2-4(\omega-1)^4}}\doteq g_{\nu}^1(\omega), \quad \text{if}\quad 0<\omega<1,\\
t_\nu(\omega)+\frac{\left[3(\omega-1)^2 + 2\nu^2\omega^4+2\nu\omega^2(\omega-1) \right]t_\nu(\omega)
-8(\omega-1)^3}{\sqrt{\left[3(\omega-1)^2 + 2\nu^2\omega^4+2\nu\omega^2(\omega-1)\right]^2-4(\omega-1)^4}}\doteq g_{\nu}^2(\omega), \quad \text{if}\quad 1<\omega<2,
\end{array}\right.
\end{equation}
where
\begin{align*}
r_\nu(\omega) &= 3\, \left( \omega -1 \right) ^{2}+2\,{\nu}^{2}{\omega}^{4}+2\,\nu{\omega
}^{2} \left( 1-\omega \right),\\
s_\nu(\omega)&=6(\omega-1)+8\nu^2\omega^3+2\nu(2\omega-3\omega^2),\\
t_\nu(\omega)&=6(\omega-1)+8\nu^2\omega^3+2\nu(3\omega^2-2\omega).
\end{align*}
It is easy to check that $g_{\nu}^2(\omega)>0$ for $1<\omega<2$ and $0<\nu<1$. Thus, for any $0<\nu<1$, $g_{\nu}(\omega)$ is strictly monotonously increasing when $1<\omega<2$. For any $0<\nu<1$, now we turn to consider $g_\nu^1(\omega)$ in \eqref{eq:dgg} with $0<\omega<1$. By direct computation, we have
\begin{align}\label{eq:g0}
g_\nu^1(0)&=-6-\frac{5\sqrt{2}}{2}<0,\\\label{eq:g1}
g_\nu^1(1)&=4v(4v-1).
\end{align}
It follows from \eqref{eq:g1} that
\begin{align}\label{eq:g1n}
g_\nu^1(1)\le 0, \quad \text{if} \quad 0<\nu\le \frac{1}{4},\\\label{eq:g1p}
g_\nu^1(1)> 0, \quad\text{if} \quad \frac{1}{4} < \nu<1.
\end{align}
In addition, we have
\begin{equation}\label{eq:dg}
 (g_\nu^1)'( \omega) = s'_\nu( \omega) + \frac{g_\nu^{1u}(\omega)}{g_\nu^{1l}(\omega)}-\frac{g_\nu^{2u}(\omega)}{ [g_\nu^{1l}(\omega)]^3},
\end{equation}
where
\begin{align*}
g_\nu^{1u}(\omega)&= [s_\nu(\omega)]^2 + s'_\nu( \omega) r_\nu(\omega) - 24(\omega - 1)^2,\\
g_\nu^{1l}(\omega)&= \sqrt{[r_\nu(\omega)]^2 - 4(\omega - 1)^4},\\
g_\nu^{2u}(\omega)&= [r_{\nu}(\omega)s_\nu(\omega)-8(\omega - 1)^3]^2.
%g_\nu^{2l}&= ([r_{\nu}(\omega)]^2 -4(\omega - 1)^4)^{\frac{3}{2}}.
\end{align*}

\begin{figure}[t]
{\centering
\begin{tabular}{ccc}
\hspace{-0.3 cm}
\resizebox*{0.50\textwidth}{0.26\textheight}{\includegraphics{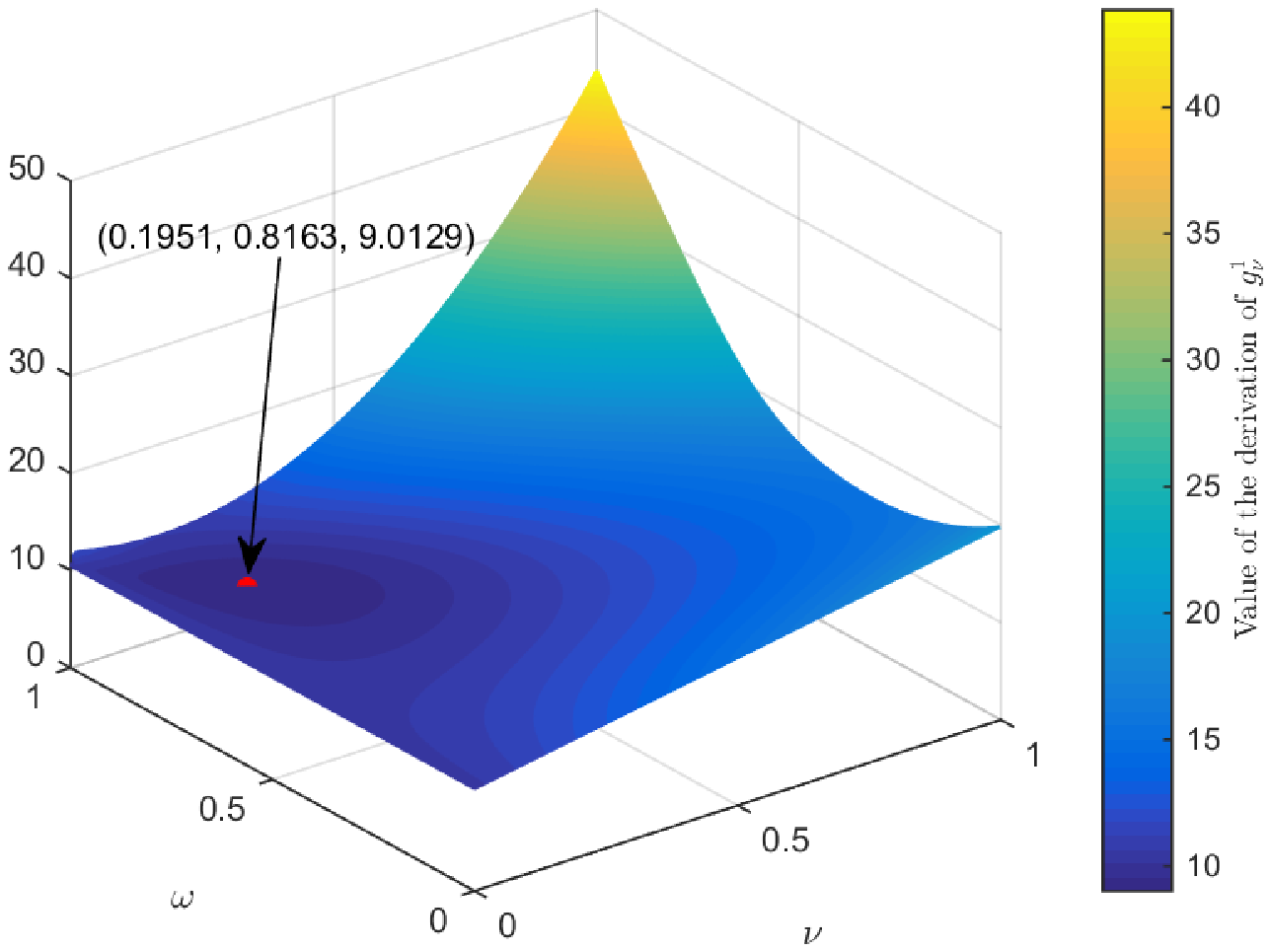}}
& & \hspace{-0.9 cm}
\resizebox*{0.50\textwidth}{0.26\textheight}{\includegraphics{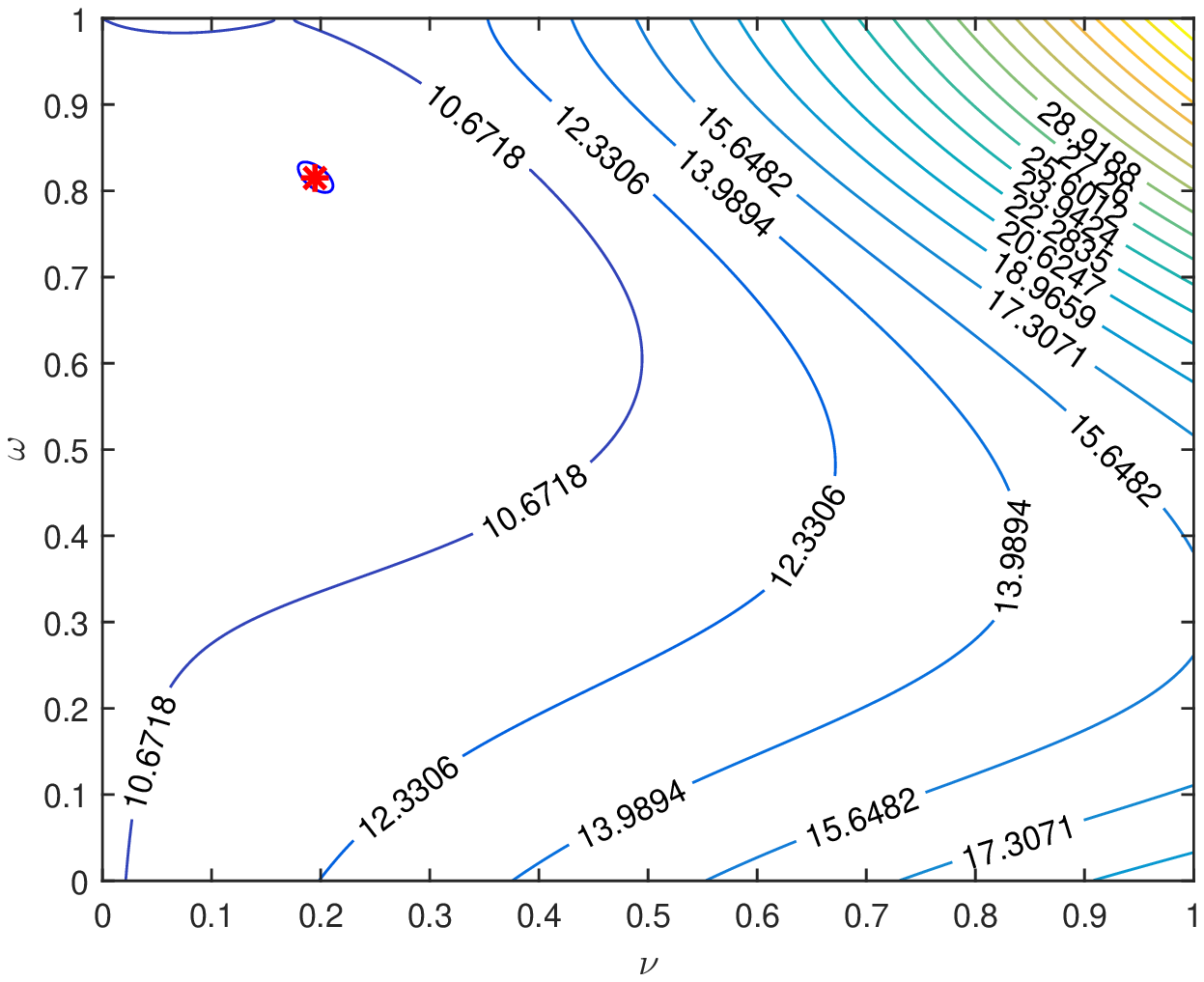}}
\end{tabular}\par
}\vspace{-0.15 cm}
\caption{Surface and contour for $(g_\nu^1)'(\omega)$ with $\nu = [0.001:0.001:0.999]$ and $\omega = [0.001:0.001:0.999]$ (The red star-dots on the surface and the contour are $(\nu,\,\omega,\,(g_\nu^1)'(\omega))=(0.1951,\,0.8163,\,9.0129)$ and $(\nu,\,\omega)=(0.1951,\,0.8163)$, respectively).}
\label{Fig:dg1}
\end{figure}

Numerically, we have always found that $(g_\nu^1)'(\omega)>0$ for $\omega\in (0,1)$ and $\nu\in (0,1)$. In fact, we turn to consider the following constrained optimization problem
\begin{equation}\label{eq:opp}
\begin{array}{rl}
\min\limits_{\omega,\nu}& (g_\nu^1)'(\omega)\\

\text{s.t.}&\left\{\begin{aligned}
&0<\omega<1,\\
&0<\nu<1.
\end{aligned}\right.
\end{array}
\end{equation}
For the constrained optimization~\eqref{eq:opp}, we simply use MATLAB's routine $\textbf{fminsearch}$ with initial value $[0.3,0.4]$ to find the solution $[\omega,\nu] \approx [0.8163, 0.1951]$ with minimum approximately equaling to $9.0129$. In Figure~\ref{Fig:dg1}, we plot the surface and contour for $(g_\nu^1)'(\omega)$ with $\nu = [0.001:0.001:0.999]$ and $\omega = [0.001:0.001:0.999]$ (MATLAB expression), which intuitively verifies the result. We should like to mention here that we have tried many ways to mathematically prove this result but without success and thus we have to left it as an open question.

\begin{figure}[t]
{\centering
\begin{tabular}{ccc}
\hspace{-0.3 cm}
\resizebox*{0.50\textwidth}{0.26\textheight}{\includegraphics{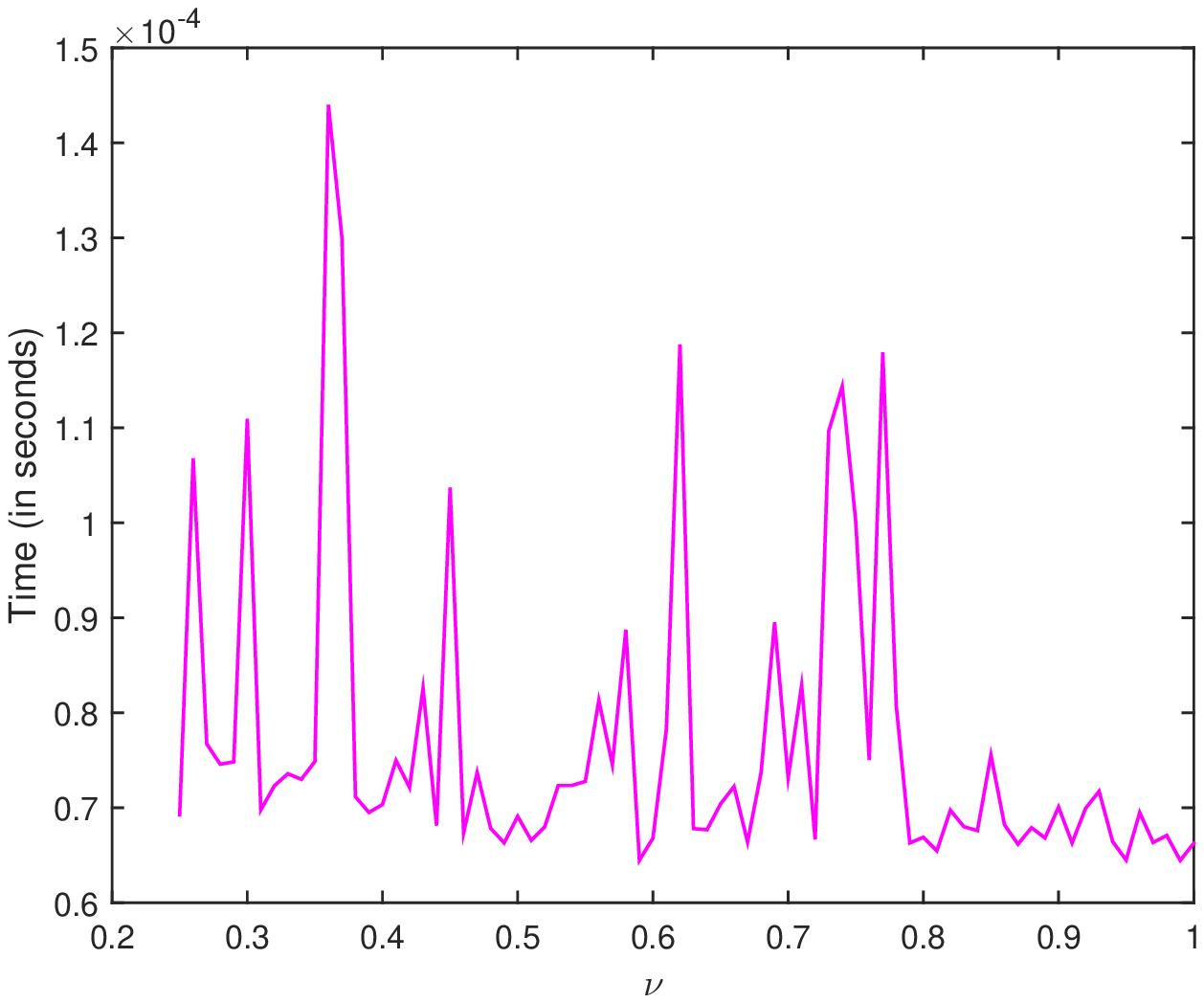}}
& & \hspace{-0.9 cm}
\resizebox*{0.50\textwidth}{0.26\textheight}{\includegraphics{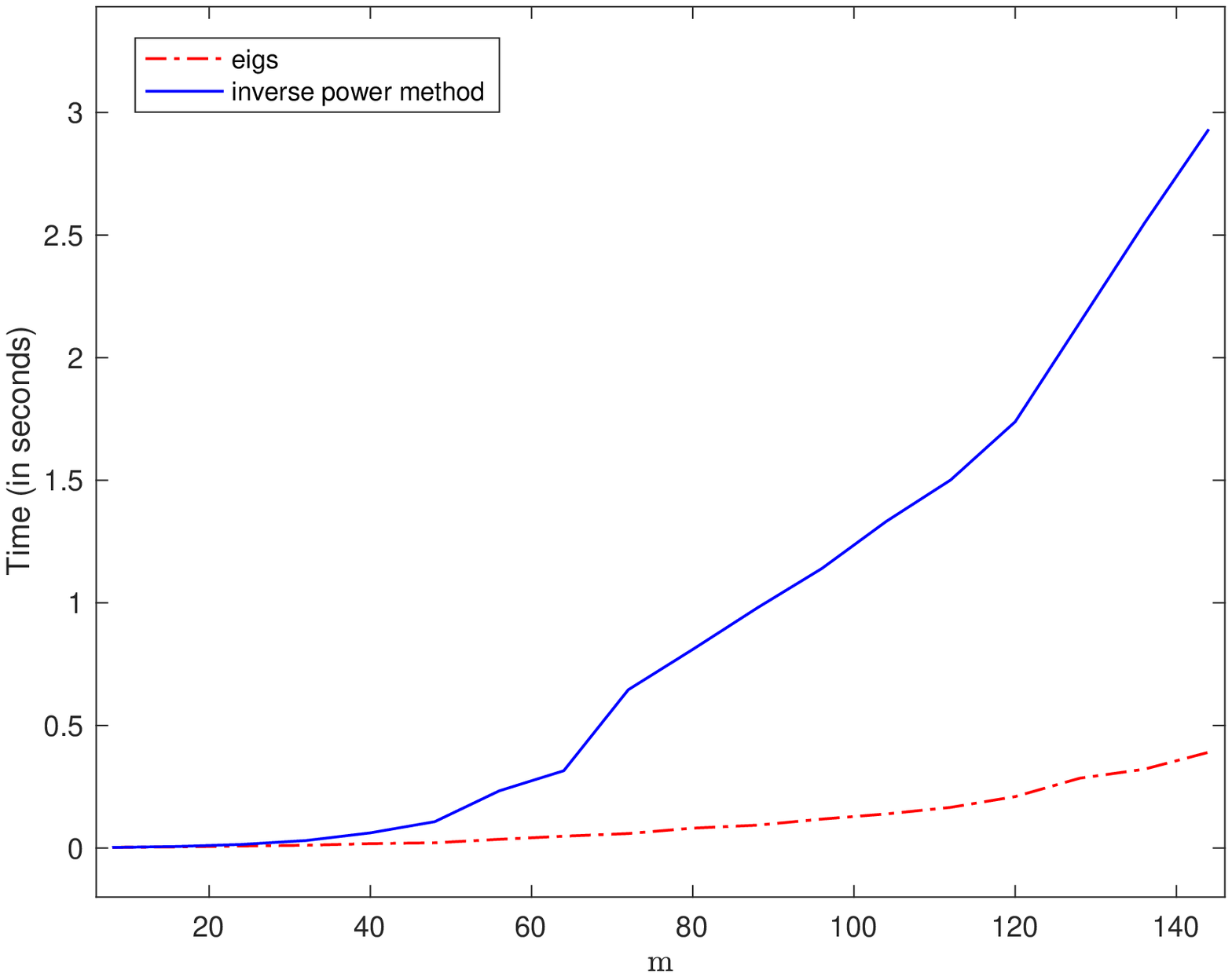}}
\end{tabular}\par
}\vspace{-0.15 cm}
\caption{Left: the elapsed CPU time for computing the root of $g_\nu^1(\omega)$ in $(0,1)$ by using bisection method (the initial ends of the interval is $b_1 = 0$ and $b_2 = 1$ and the termination criterion is $|g_\nu^1(\omega)|\le 10^{-12}$ or the updated ends of the interval $b_2-b_1\le 10^{-12}$); Right: the elapsed CPU time for computing $\nu =\|A^{-1}\|_2$ with matrices in Example~\ref{Example4.2} (see Section~\ref{sec:numer} for detail).}
\label{Fig:time4opt}
\end{figure}

\begin{figure}[t]
{\centering
\begin{tabular}{ccc}
\hspace{-0.3 cm}
\resizebox*{0.50\textwidth}{0.26\textheight}{\includegraphics{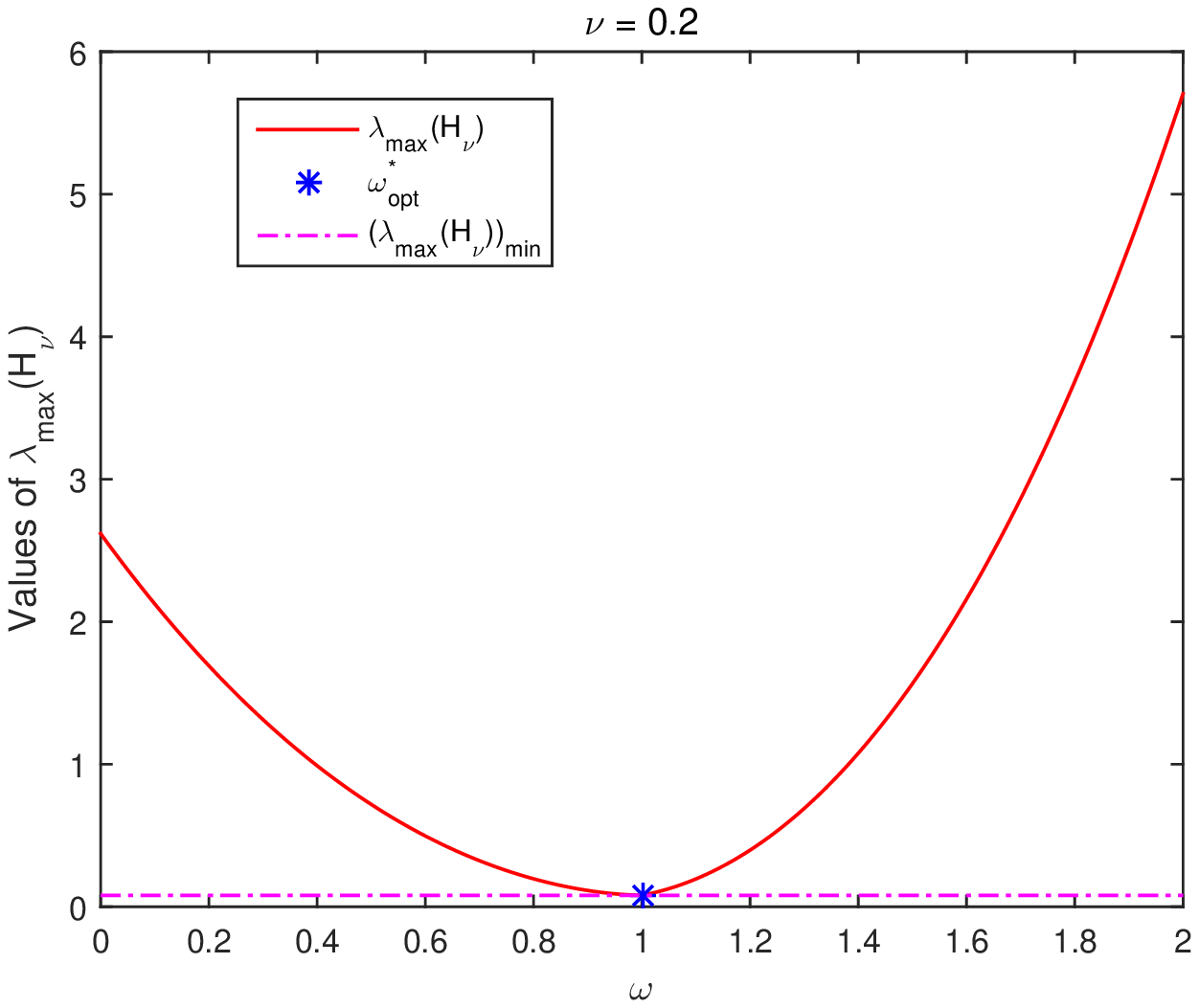}}
& & \hspace{-0.9 cm}
\resizebox*{0.50\textwidth}{0.26\textheight}{\includegraphics{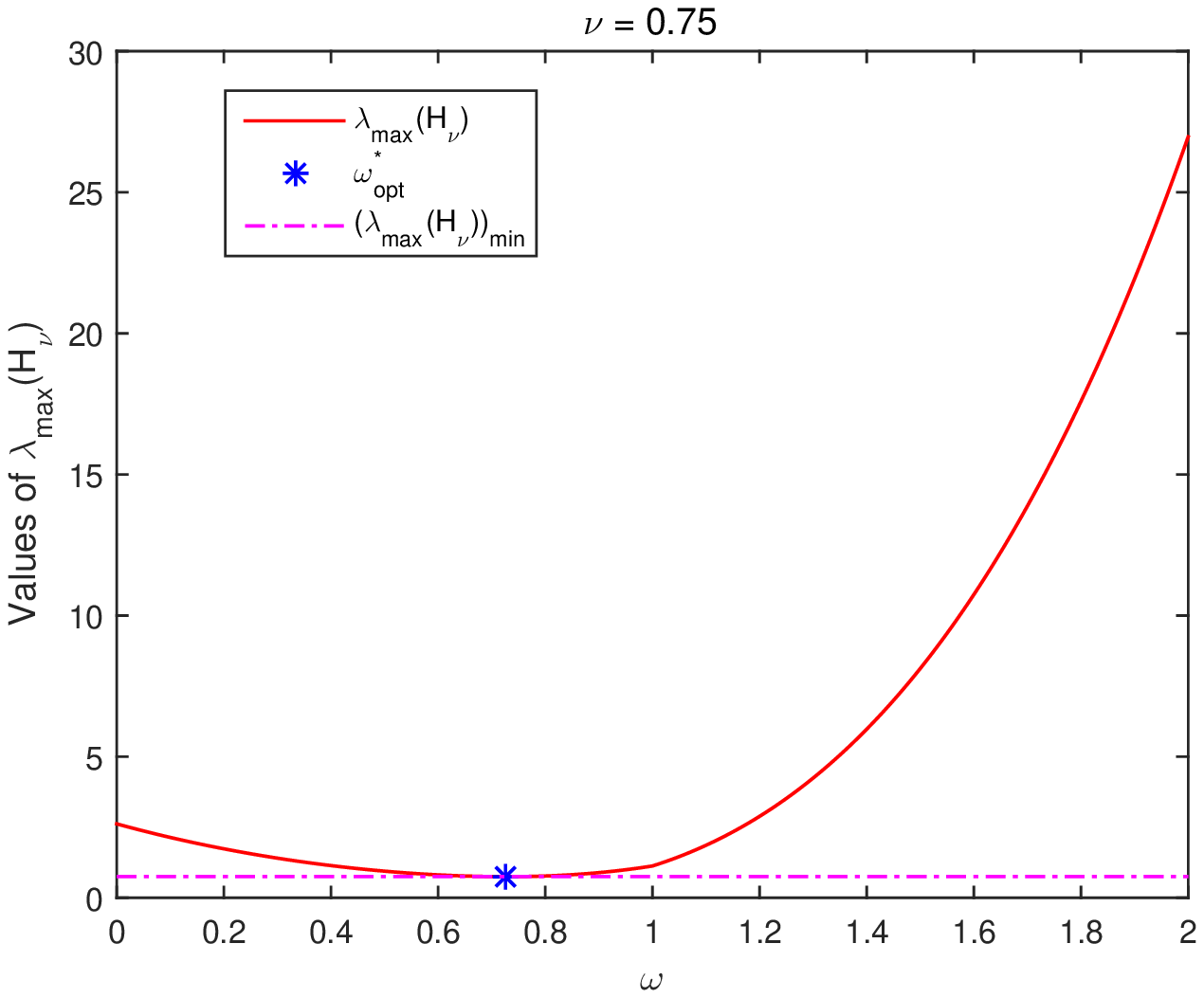}}
\end{tabular}\par
}\vspace{-0.15 cm}
\caption{Curves of $\lambda_{\max}(H_\nu(\omega))$.}
\label{Fig:lambda}
\end{figure}

Based on this observation, we assume that $(g_\nu^1)'(\omega) > 0$ for $\omega\in (0,1)$ and $\nu\in (0,1)$ in the following discussion. Combining it (which means that $g_\nu^1(\omega)$ is strictly monotonously increasing in $\omega\in (0,1)$) with \eqref{eq:g0}, \eqref{eq:g1n} and \eqref{eq:g1p}, we can conclude that $g_\nu^1(\omega)$ has at most one root when $0<\omega<1$. If the root exists (when $\frac{1}{4} < \nu<1$), which is denoted by $\omega_{opt}$, then $g_{\nu}$ is strictly monotonously decreasing in $(0, \omega_{opt})$ and strictly monotonously increasing in $(\omega_{opt},2)$ and the optimal iteration parameter which minimizes $\|T_\nu(\omega)\|_2$ is $\omega^*_{opt}=\omega_{opt}$. Otherwise, if $0<\nu\le\frac{1}{4}$, then $g_{\nu}$ is strictly monotonously decreasing in $(0, 1)$ and strictly monotonously increasing in $(1,2)$ and the optimal iteration parameter which minimizes $\|T_\nu(\omega)\|_2$ is $\omega^*_{opt}=1$. In conclusion, the optimal iteration parameter is
\begin{equation}\label{eq:opt}
\omega^*_{opt}=\left\{\begin{array}{ll}
\omega_{opt},\quad \text{if}\quad \frac{1}{4} < \nu<1,\\
1,\quad \text{if}\quad 0<\nu\le\frac{1}{4}.
\end{array}\right.
\end{equation}
The root of $g_\nu^1(\omega)$ located in $(0,1)$, $\omega_{opt}$, can be efficiently numerically calculated by the classical bisection method and the elapsed CPU time (the average of ten tests) is shown in Figure~\ref{Fig:time4opt} (the left plot), from which we find that the root can be obtained within $1.5\times 10^{-4}$ seconds provided that $\nu$ is known. Figure~\ref{Fig:lambda} shows some curves of $\lambda_{\max}(H_\nu(\omega))$ which visually demonstrates our claims mentioned above.

\begin{figure}[t]
{\centering
\begin{tabular}{ccc}
\hspace{-0.3 cm}
\resizebox*{0.50\textwidth}{0.26\textheight}{\includegraphics{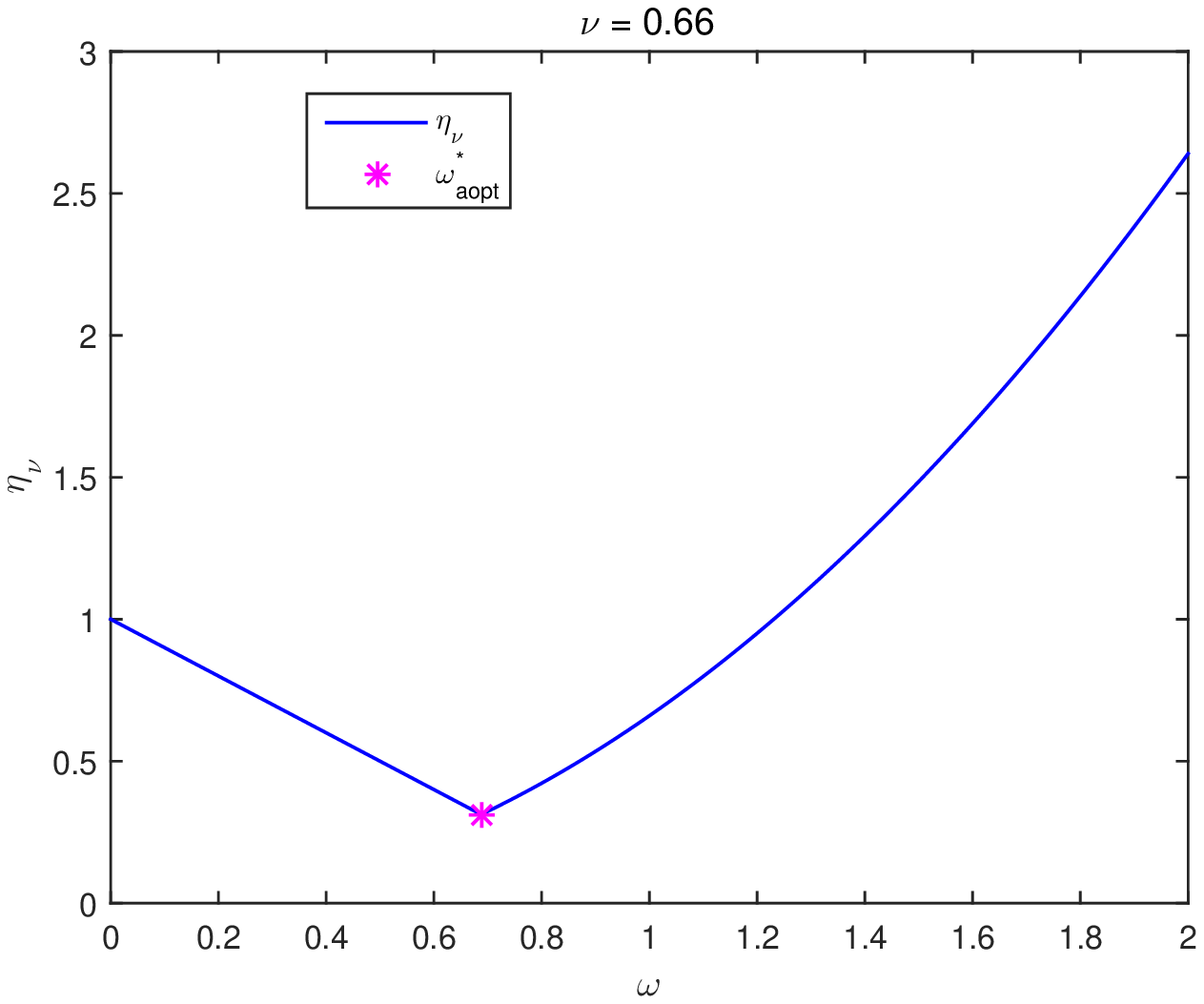}}
& & \hspace{-0.9 cm}
\resizebox*{0.50\textwidth}{0.26\textheight}{\includegraphics{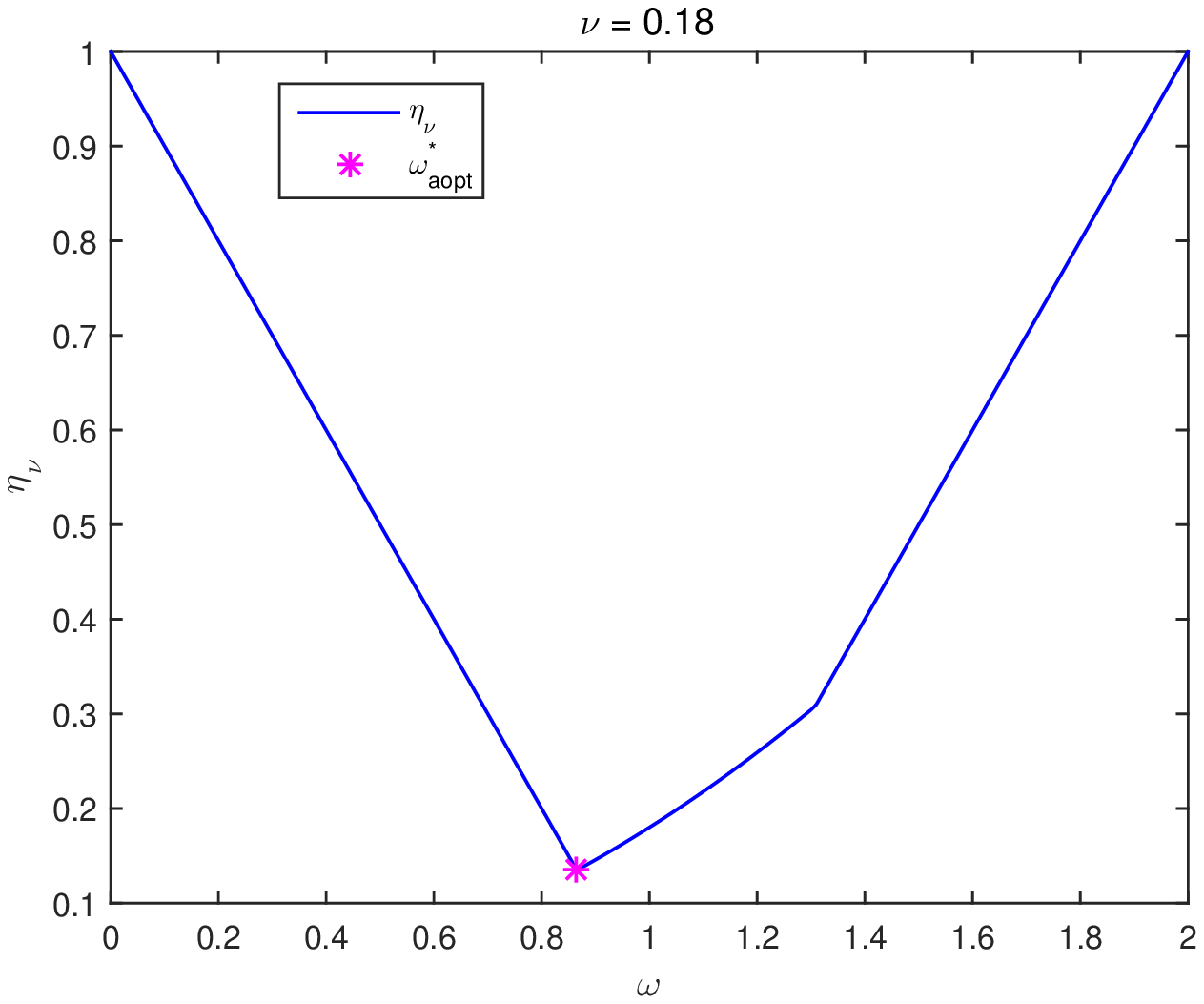}}
\end{tabular}\par
}\vspace{-0.15 cm}
\caption{Curves of $\eta_\nu(\omega)$.}
\label{Fig:eta}
\end{figure}

\begin{figure}[t]
{\centering
\begin{tabular}{ccc}
\hspace{-0.3 cm}
\resizebox*{0.50\textwidth}{0.26\textheight}{\includegraphics{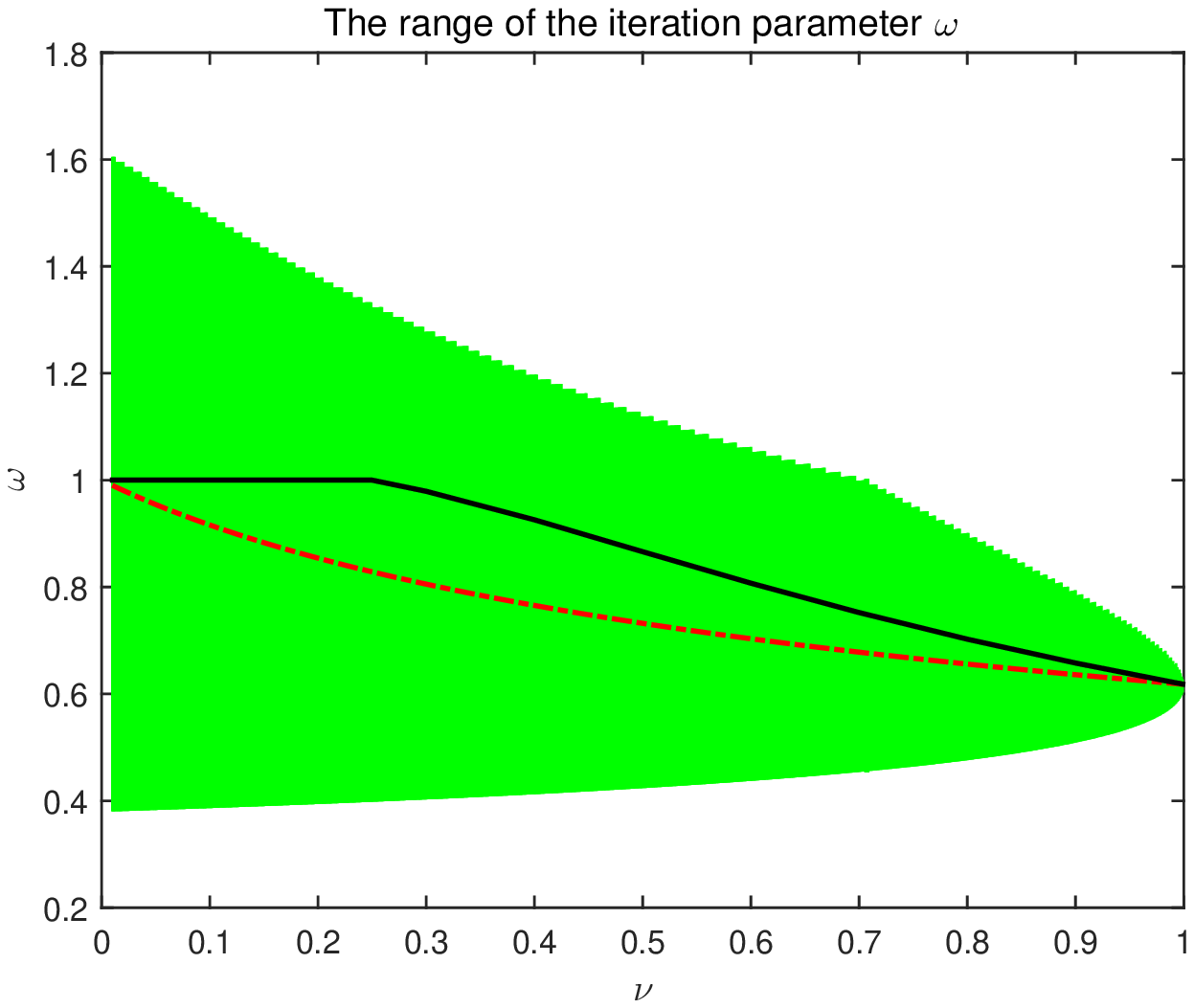}}
& & \hspace{-0.9 cm}
\resizebox*{0.50\textwidth}{0.26\textheight}{\includegraphics{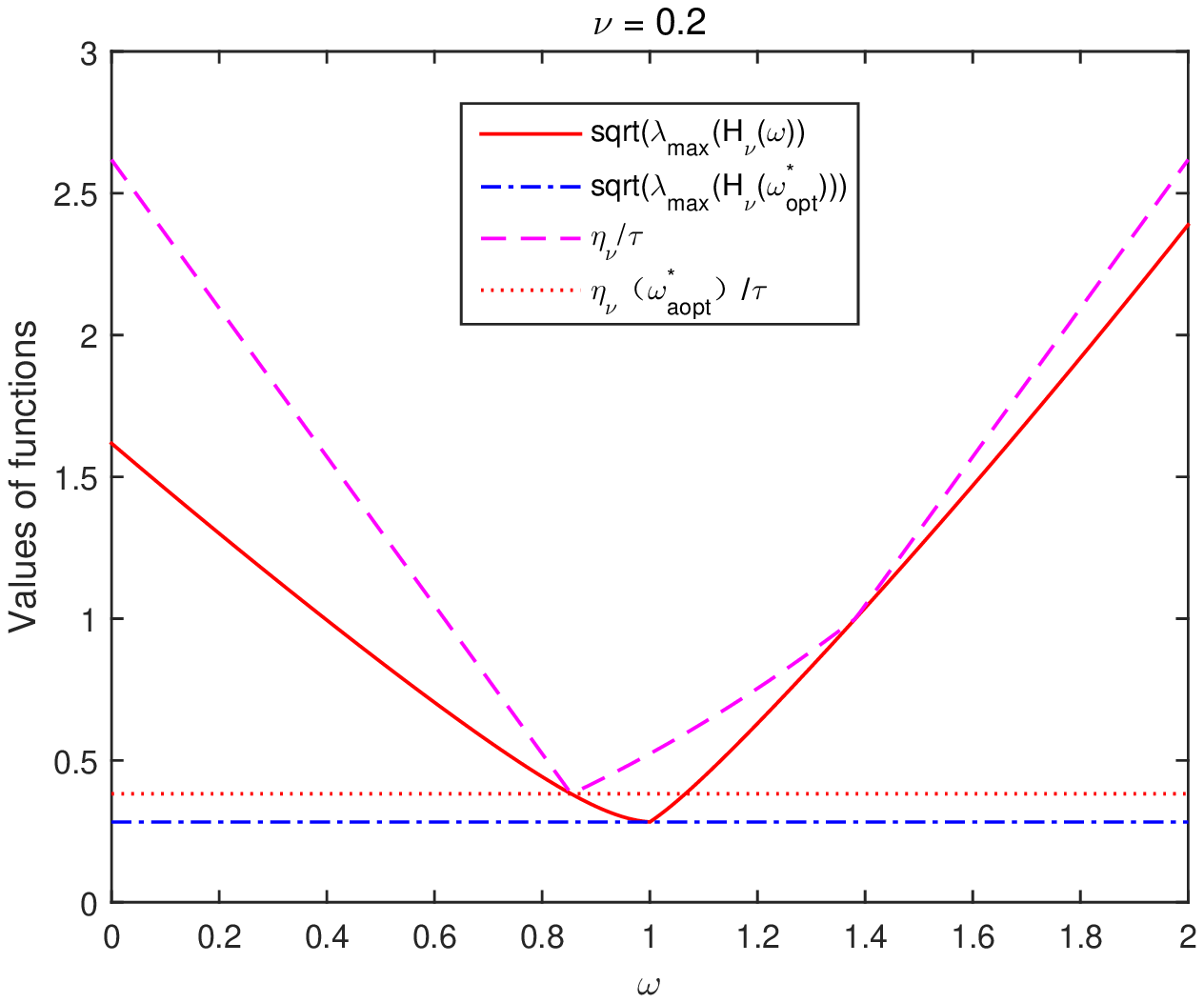}}\vspace{2ex}\\
\hspace{-0.3 cm}
\resizebox*{0.50\textwidth}{0.26\textheight}{\includegraphics{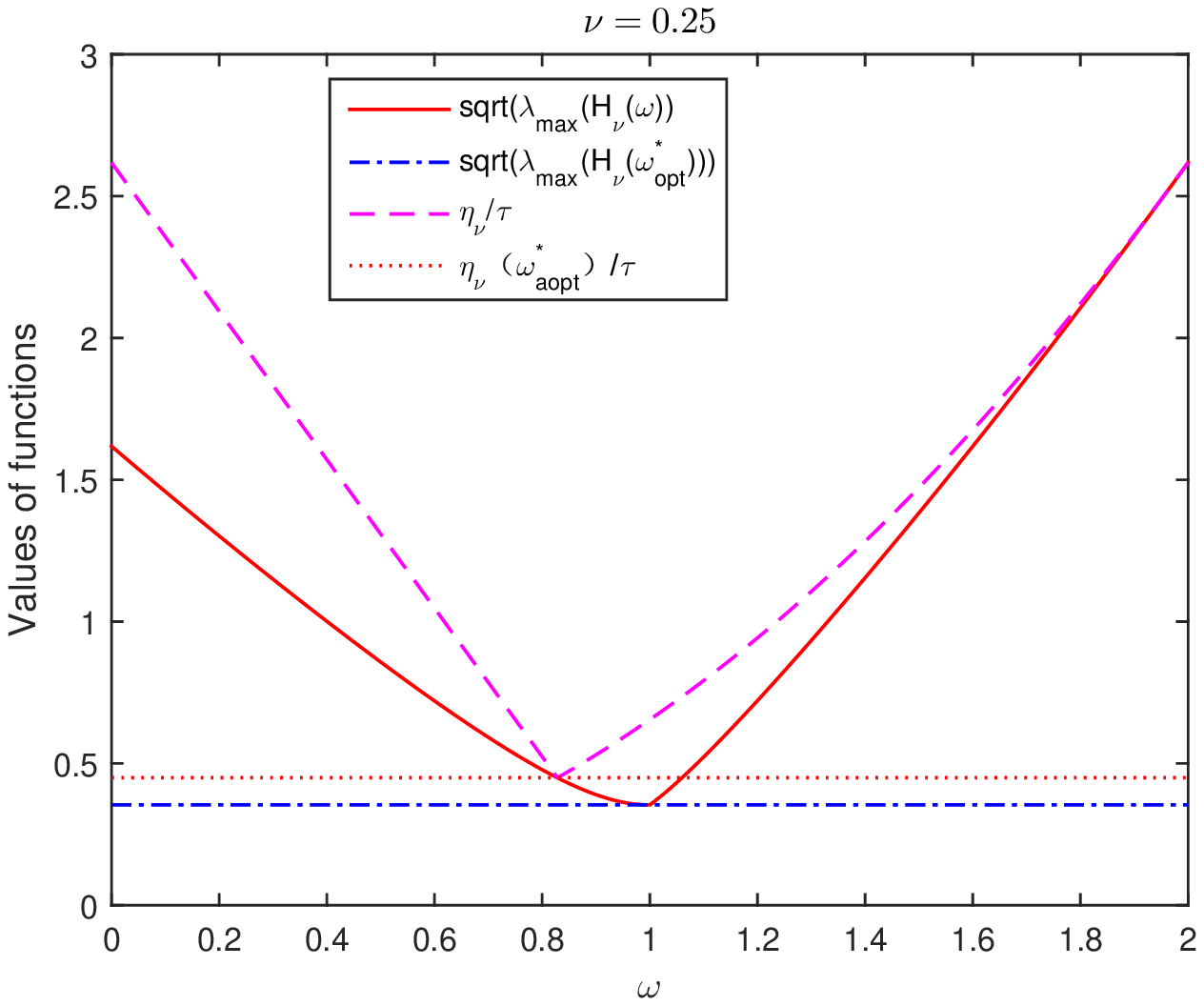}}
& & \hspace{-0.9 cm}
\resizebox*{0.50\textwidth}{0.26\textheight}{\includegraphics{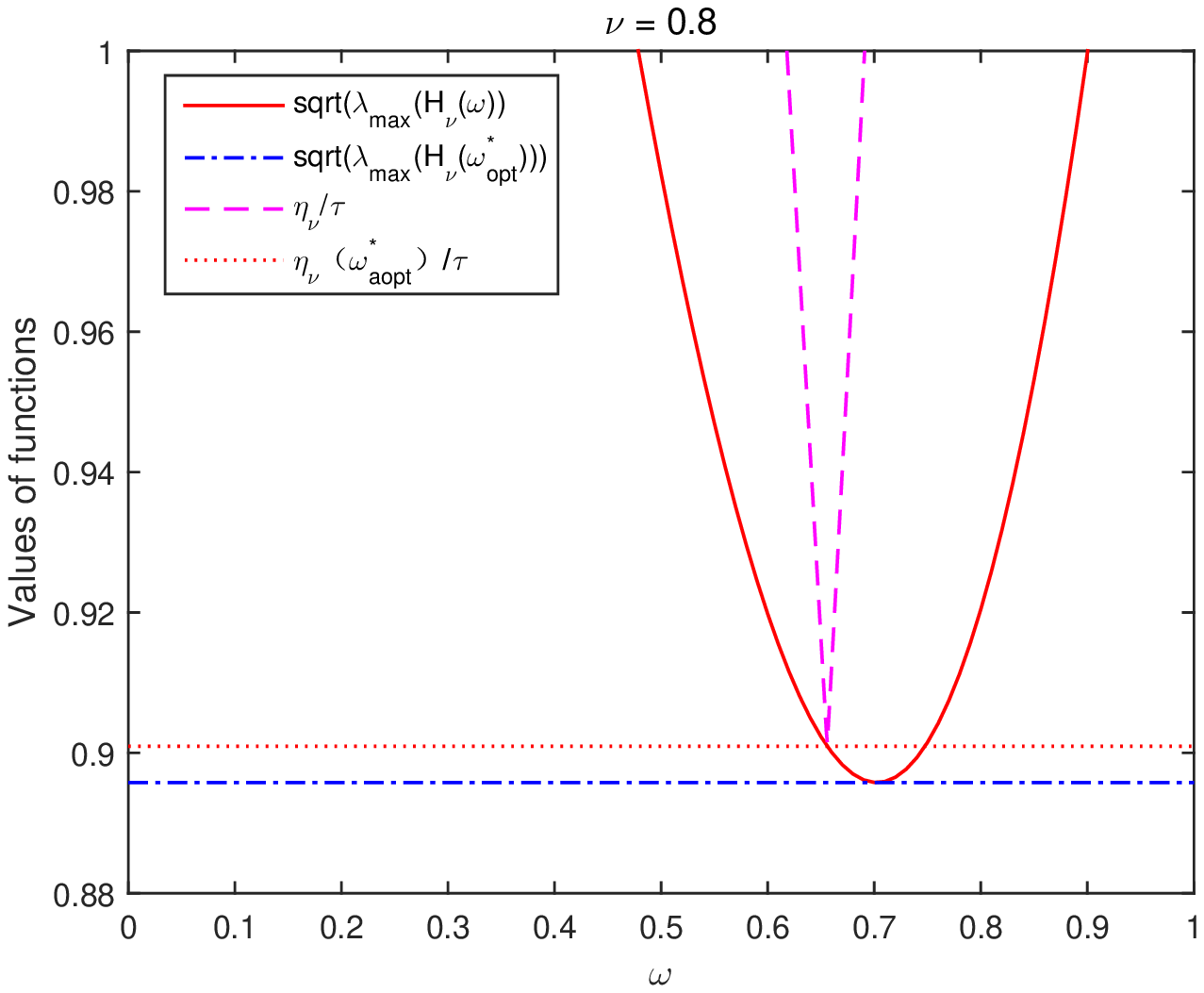}}
\end{tabular}\par
}\vspace{-0.15 cm}
\caption{Curves of optimal (top left: black solid line) and approximate optimal (top left: red dots line) iteration parameters and the comparisons of minimal points and minimums (others).}
\label{Fig:lambdavseta}
\end{figure}

Now we devote to discussing the approximate optimal iteration parameter $\omega^*_{aopt}$ which will minimize $\eta_\nu(\omega)=\max \{|1-\omega|\doteq h(\omega),\omega^2 \nu\doteq h_{\nu}(\omega) \}$ with $\omega \in (0,2)$. This $\eta_\nu(\omega)$ appears in the proof of Corollary~3.1 in \cite{kema2017} and satisfies $\|T_\nu(\omega)\|_2\le \frac{\eta_\nu(\omega)}{\tau}$. That is, $\frac{\eta_\nu(\omega)}{\tau}$ is an upper bound of $\|T_\nu(\omega)\|_2$ with $\omega \in (0,2)$, which is the reason that we call $\omega^*_{aopt}$ the approximate optimal iteration parameter. Since $h(\omega)$ is strictly monotonously decreasing in $(0,1)$ and strictly monotonously increasing in $(1,2)$ and $h_{\nu}(\omega)$ is strictly monotonously increasing in $(0,2)$ with $\nu\in (0,1)$, $\omega^*_{aopt}$ must satisfy $1-\omega=\nu \omega^2(0<\omega<2)$, that is
\begin{equation}\label{eq:waopt}
\omega^*_{aopt}(\nu)=\frac{\sqrt{4\nu+1}-1}{2\nu}.
\end{equation}
Figure~\ref{Fig:eta} intuitively confirms this result. Furthermore, from~\eqref{eq:waopt}, the approximate optimal iteration parameter $\omega^*_{aopt}(\nu)$ is strictly monotonously decreasing with respect to $\nu\in(0,1)$ and  it approaches to $1$ (the optimal one when $0<\nu\le\frac{1}{4}$) as $\nu\rightarrow 0^+$ and closes to $\frac{\sqrt{5}-1}{2}$ as $\nu\rightarrow 1^-$. Thus, $\omega^*_{aopt}(\nu)$ must located in the range of $\omega$ determined by conditions~\eqref{eq:ncond1}-\eqref{eq:ncond3} for $\nu\in(0,1)$. On the other hand, since the optimal iteration parameter $\omega^*_{opt}$ must satisfy the conditions~\eqref{eq:cond1}, it also must located in the range of $\omega$ determined by conditions~\eqref{eq:ncond1}-\eqref{eq:ncond3} and it closes to $\frac{\sqrt{5}-1}{2}$ as $\nu\rightarrow 1^-$. Furthermore, from the definition of the approximate optimal iteration parameter, the value of $\sqrt{\lambda_{\max}(H_\nu(\omega))}$ is no larger than that of $\frac{\eta_\nu(\omega)}{\tau}$ and thus $\sqrt{\lambda_{\max}(H_\nu(\omega^*_{opt}))}\le \frac{\eta_\nu(\omega^*_{aopt})}{\tau}$, which implies that the SOR-like iteration method with optimal iteration parameter \eqref{eq:opt} will converge no slower than that with the approximate optimal iteration parameter~\eqref{eq:waopt}. Indeed, when $\nu\in (0,\frac{1}{4}]$, we have $\omega^*_{aopt}<1=\omega^*_{opt}$, combining this with the strictly monotonous properties of $\frac{\eta_\nu(\omega)}{\tau}$ and $\sqrt{\lambda_{\max}(H_\nu(\omega))}$, we can conclude that $\sqrt{\lambda_{\max}(H_\nu(\omega^*_{opt}))}< \frac{\eta_\nu(\omega^*_{aopt})}{\tau}$. When $\nu\in (\frac{1}{4},1)$, the analytical formulation of $\omega^*_{opt}$ is unknown, and,  numerically, we also find that $\omega^*_{aopt}<\omega^*_{opt}$ whenever $\omega^*_{opt}$ is obtained. Thus, $\sqrt{\lambda_{\max}(H_\nu(\omega^*_{opt}))}< \frac{\eta_\nu(\omega^*_{aopt})}{\tau}$ for $\nu \in (0,1)$, which means that the SOR-like iteration method with optimal iteration parameter \eqref{eq:opt} converges faster than that with the approximate optimal iteration parameter~\eqref{eq:waopt}. Figure~\ref{Fig:lambdavseta} intuitively clarifies most of our arguments here. Moreover, our numerical results in the next section will further demonstrate these arguments.

%%%%%%%%%%%%%%%%%%%%%%%%%%%%%%%%%%

\begin{rem}
In \cite{guwl2019}, the authors give some  convergence conditions from the involved iteration matrix of the SOR-like iteration method for solving the AVE~\eqref{eq:ave}, which are different from the results of \cite{kema2017}. In addition, if all the eigenvalues of $D(x^{(k+1)})A^{-1}$ are positive, then the optimal iteration parameter which minimizes the spectral radius of the iteration matrix is $\omega^*_o=\frac{2}{1+\sqrt{1-\varrho}}$ with $\varrho=\rho(D(x^{(k+1)})A^{-1})$~\cite{guwl2019}. Apparently, the optimal iteration parameter is iteration-dependent and it will be expensive to compute it in every iterative step, and the authors use the approximate optimal one $\omega_o=\frac{2}{1+\sqrt{1-\tilde{\varrho}}}$ with $\tilde{\varrho}=\rho(A^{-1})$ in their numerical experiments. Though the approximate optimal parameter is iteration-independent, however, as it will be shown in the next section, $\omega_o$ will be out of the range of $\omega$ determined by \eqref{eq:ncond1}-\eqref{eq:ncond3}, which means that the SOR-like method with this $\omega_o$ will be divergent. Finally, we want to point out that the SOR-like iteration method may not be a stationary iterative method $($it will be if all the iteration $x^{(k)}$ possess the same sign for $k = 1,2,\cdots$$)$. Thus, the spectral radius of the iteration matrix less than $1$ will not guarantee the convergence, a counterexample can be found in \cite[Pages~6--7]{gual2015}.
\end{rem}

\begin{rem}\label{rem:opt}
From our analysis above, $\omega_{opt}^*$ is the optimal iteration parameter which minimizes $\|T_\nu(\omega)\|_2$, an upper bound of the linear convergence factor for the SOR-like iteration method in the metric $\| |\cdot |\|_{\omega}$ of $(e_{k}^x,e_{k}^y)$ $($see \eqref{ie:cov}$)$. However, it may not be the really optimal one since the upper bound may not be tight; see the numerical results of Example~\ref{Example4.3} for an instance.
\end{rem}
%%%%%%%%%%%%%%%%%%%数值实验%%%%%%%%%%%%%%%%%%%%%%%%%%
\section{Numerical experiments}\label{sec:numer}
\qquad In this section, we will present four numerical examples to illustrate the superior performance of the SOR-like iteration method with our optimal iteration parameter for solving several AVEs~\eqref{eq:ave}. Five algorithms will be tested.
\begin{enumerate}

  \item SORLo: the SOR-like iteration method with the approximate optimal iteration parameter $\omega_{o}=\frac{2}{1+\sqrt{1-\tilde{\varrho}}}$  proposed in \cite{guwl2019}. Here, $\tilde{\varrho} = \rho(A^{-1})$.

  \item SORLaopt: the SOR-like iteration method with the approximate optimal iteration parameter $\omega^*_{aopt}$ defined as in \eqref{eq:waopt}.

  \item SORLopt: the SOR-like iteration method with the optimal iteration parameter $\omega^*_{opt}$ defined as in \eqref{eq:opt}.

  \item SORLno: the SOR-like iteration method with the numerically optimal iteration parameter $\omega_{no}$, which is selected from $\omega=[0.001:0.001:1.999]$ and is the first one to reach the minimal number of iteration of the method. Practically, determining optimal iteration parameter in this way is unwise because of that it is always very time consuming, especially for large scale problems (see the following numerical results for detail). However, it is still necessary to be tested here in order to demonstrate the claim in Remark~\ref{rem:opt}.

  \item NT: the generalized Newton method \cite{mang2009}. Given an
            initial guess $x^{(0)}\in \mathbb{R}^n$, the iteration sequence $\{x^{(k)}\}$ is generated by
           \begin{equation}\label{eq:newton}
               \left[A-D(x^{(k)})\right]x^{(k+1)} = b
            \end{equation}
             for $k = 0,\,1,\,\cdots$ until convergence.

\end{enumerate}

Note that in all these algorithms, the main task per iteration is solving a system of linear equations. The main difference between NT and SOR-like algorithms is that the coefficient matrices: In NT, the matrix is $A-D(x^{(k)})$, which is varying along with the iteration; while in the SOR-like methods, the matrix is $A$, a fixed one. Solving systems $Ax=d^k$ for different $k$ only requires that $A$ be factored once, at a cost of $O(n^3)$ flops for the worst case ($A$ is dense); while per iteration, its cost reduces to $O(n^2)$ flops (matrix-vector products). Hence, if there are $K$ iteration steps, the cost is $O(n^3)+KO(n^2)$. However, for the NT, the cost will be $K O(n^3)$ for the worst case. Thus, if $n$ is large and $K$ is far smaller than $n$, the SOR-like methods will be competitive in terms of CPU time. In this paper, the tested methods are implemented in conjunction with Cholesky factorization (for the former three examples) or LU factorization (for the last example). Specifically, we use
 $
 dA ={ \bf decomposition}(A,\text{`chol'})~\text{or}~dA ={ \bf decomposition}(A,\text{`lu'})
 $
 to generate the Cholesky or LU decomposition of $A$, where ${\bf decomposition}$ is the routine in MATLAB, which returns the corresponding decomposition of a matrix $A$ that can be used to solve the linear system $Ax=b$ efficiently. The call $x=dA\backslash b$ returns the same vector as $A\backslash b$, but is typically faster.

We also should note that the SOR-like method is one parameter-dependent and computing the optimal or approximate optimal parameters involved $\|A^{-1}\|_2$ or $\rho(A^{-1})$, and thus it may be time consuming, especially for large sparse problems. For symmetric positive definite matrix $A$, we have $\nu = \tilde{\varrho} = \frac{1}{\lambda_{\min}(A)}$. In the literature, inverse power method is available for computing $\lambda_{\min}(A)$; see~\cite{bai2000s} for instance. The right plot in Figure~\ref{Fig:time4opt} compares the time performance of the inverse power method (the initial point is $\textbf{1}_n$, the tolerance is $10^{-8}$ and the maximal iteration number is $2000$) and the MATLAB's routine $\textbf{eigs}(A,1,\text{`SA'})$ for computing $\|A^{-1}\|_2$ or $\rho(A^{-1})$ of the symmetric positive matrices in the following Example~\ref{Example4.2}, from which we find that the elapsed CPU time is increasing as $n$ becoming larger and, in terms of CPU time, $\textbf{eigs}$ is better than the inverse power method for this example, especially when $n$ is large. However, sometimes $\textbf{eigs}$ can't work (for instance, for Example~\ref{Example4.1}) and, in this case, we use the inverse power method in our tests.

In the numerical results, we will report ``IT" (the number of iteration), ``CPU" (the elapsed CPU time in seconds) and ``RES" (the residual error). RES is defined by
$$\text{RES} = \|Ax^{(k)}-|x^{(k)}|-b\|.$$
All tests are started from the initial zero vector (except Example~\ref{Example4.4}) and terminated if the current iteration satisfies $\text{RES} \le 10^{-8}$ (except the last problem in Example~\ref{Example4.3} and Example~\ref{Example4.4}) or the number of prescribed maximal iteration steps $k_{\text{max}}=100$ is exceeded (denoted by ``--"). In order to obtain more accurate CPU time, we run all test problems in each method five times and take the average. All computations are done in MATLAB R2017b with a machine precision $2.22\times 10^{-16}$ on a personal computer with 2.60GHz central processing unit (Intel Core i7), 16GB memory and MacOS operating system.

\begin{exam}[\!\!\cite{guwl2019}]\label{Example4.1}  Consider the AVE \eqref{eq:ave} with
$$ A = \text{tridiag}(-1,8,-1)=\left[\begin{array}{cccccc}
                                     8&-1&0&\cdots&0&0\\
                                     -1&8&-1&\cdots&0&0\\
                                     0&-1&8&\cdots&0&0\\
                                     \vdots&\vdots&\vdots&\ddots&\vdots&\vdots\\
                                     0&0&0&\cdots&8&-1\\
                                     0&0&0&\cdots&-1&8\end{array}\right]\in \mathbb{R}^{n\times n}\quad\text{and}\quad b=Ax^*-|x^*|,$$
where $x^*=[-1,1,-1,1,\cdots,-1,1]^T\in \mathbb{R}^n.$

For this example, we always have $\nu=0.1667<0.25,\,\omega_{o}=1.0455,\,\omega^*_{aopt}=0.8730,\,\omega^*_{opt}=1$ and $\omega_{no}= (0.994,0.997,0.983,0.984,0.985)$ for $n = (1000,2000,3000,4000,5000)$, respectively. Obviously, $\omega^*_{aopt}<\omega^*_{opt}$ and $0.2357=\sqrt{\lambda_{\max}(H_\nu(\omega^*_{opt}))}<0.3326=\frac{\eta(\omega^*_{aopt})}{\tau}<1$. In addition, the range of the iteration parameter $\omega$ is $\Omega_1=(0.3938,1.4184)$ and $\omega_{o},\,\omega_{no},\, \omega^*_{aopt},\, \omega^*_{opt} \in \Omega_1$.

Tabel~\ref{table1} reports the numerical results for Example~\ref{Example4.1}, from which we clearly find that SORLopt converges faster than both SORLaopt and SORLo. SORLopt and SORLno take the same number of iteration for this example, but SORLopt is better than SORLno in terms of CPU time and the reason is that determining the numerically optimal iteration parameter needs more CPU than computing the optimal iteration parameter. SORLo is better than SORLaopt with respect to IT and CPU. In regard to IT, CPU and RES, NT is the best one among the tested methods for this example.
\setlength{\tabcolsep}{1.5pt}
\begin{table}[!h]\small
\centering
\caption{Numerical results for Example~\ref{Example4.1}.}\label{table1}
\begin{tabular}{|c|c|ccccc|}\hline
 \multirow{2}{*}{Method}  &  & \multicolumn{5}{|c|}{$n$}  \\ \cline{3-7}
  &     &  $1000$ & $2000$ & $3000$ & $4000$ & $5000$ \vextra\\ \hline\hline
%%%%%%%%%%%%%%%%%%%%%%%%%%%%%%%%%%%%%%%%%%%%%%%%%%%%%%%%%%%%%%%%%%%%%%%%%%%%%%%%%%
\multirow{3}{*}{SORLo} & IT  & $16$ & $16$ & $17$ & $17$& $17$ \\
                 & CPU & $0.0465$ & $0.0733$ & $0.1070$ & $0.1522$& $0.1942$ \\
           & RES & $6.7263\times 10^{-9}$ & $9.5154\times 10^{-9}$  & $2.5980\times 10^{-9}$ &  $3.0001\times 10^{-9}$& $3.3543\times 10^{-9}$\\\hline
 \multirow{3}{*}{SORLaopt} & IT  & $20$ & $20$ & $20$ & $20$& $20$\\
                 & CPU & $0.0467$ & $0.0737$ & $0.1084$ & $0.1576$& $0.1985$\\
            & RES &$4.1299\times 10^{-9}$ & $5.8421\times 10^{-9}$ & $7.1557\times 10^{-9}$ & $8.2630\times 10^{-9}$& $9.2385\times 10^{-9}$\\\hline
  \multirow{3}{*}{SORLopt} & IT  & $12$ & $12$ & $13$ & $13$& $13$\\
                 & CPU & $0.0446$ & $0.0729$ & $0.1062$ & $0.1490$& $0.1921$\\
            & RES &$6.8073\times 10^{-9}$ &$9.6301\times 10^{-9}$ & $1.1900\times 10^{-9}$ & $1.3741\times 10^{-9}$& $1.5364\times 10^{-9}$ \\\hline
     \multirow{3}{*}{SORLno} & IT  & $12$ & $12$ & $13$ & $13$& $13$\\
   & CPU & $4.1239$ & $6.9784$ & $9.8242$ & $12.6611$& $16.5506$\\
   & RES &$9.3737\times 10^{-9}$ &$9.8598\times 10^{-9}$ & $9.6529\times 10^{-9}$ & $9.8334\times 10^{-9}$& $9.6793\times 10^{-9}$ \\\hline
     \multirow{3}{*}{NT} & IT  & $\textbf{2}$ & $\textbf{2}$ & $\textbf{2}$ & $\textbf{2}$& $\textbf{2}$\\
   & CPU & $\textbf{0.0011}$ & $\textbf{0.0016}$ & $\textbf{0.0020}$ & $\textbf{0.0025}$& $\textbf{0.0030}$\\
   & RES &$\boldsymbol{4.6998\times 10^{-15}}$ &$\boldsymbol{4.6998\times 10^{-15}}$ & $\boldsymbol{4.6998\times 10^{-15}}$ & $\boldsymbol{4.6998\times 10^{-15}}$& $\boldsymbol{4.6998\times 10^{-15}}$ \\ \hline
\end{tabular}
\end{table}

\end{exam}
%%%%%%%%%%%%%%%%%%%%%%%%%%%%%%%%%%%%%%%%%%%%%%%%%%%%%%%%

\begin{exam}[\!\!\cite{guwl2019}]\label{Example4.2} Consider the AVE \eqref{eq:ave} with
$$A =\text{Tridiag}(-I_m,S_m,-I_m)=\left[\begin{array}{cccccc}
                                     S_m&-I_m&0&\cdots&0&0\\
                                     -I_m&S_m&-I_m&\cdots&0&0\\
                                     0&-I_m&S_m&\cdots&0&0\\
                                     \vdots&\vdots&\vdots&\ddots&\vdots&\vdots\\
                                     0&0&0&\cdots&S_m&-I_m\\
                                     0&0&0&\cdots&-I_m&S_m\end{array}\right]\in \mathbb{R}^{n\times n},$$
$$S_m = \text{tridiag}(-1,8,-1)=\left[\begin{array}{cccccc}
                                     8&-1&0&\cdots&0&0\\
                                     -1&8&-1&\cdots&0&0\\
                                     0&-1&8&\cdots&0&0\\
                                     \vdots&\vdots&\vdots&\ddots&\vdots&\vdots\\
                                     0&0&0&\cdots&8&-1\\
                                     0&0&0&\cdots&-1&8\end{array}\right]\in \mathbb{R}^{m\times m}$$
and $b=Ax^*-|x^*|$, where $x^*=[-1,1,-1,1,\cdots,-1,1]^T\in \mathbb{R}^n$. Here, we have $n=m^2$.

The parameters for this example are displayed in Table~\ref{table-range}, from which we can find that $\omega^*_{aopt}<\omega^*_{opt}$, $\omega_{o},\,\omega_{no},\, \omega^*_{aopt},\,\omega^*_{opt} \in \Omega_1$ and the larger the value of $\nu$ is, the smaller the range of $\omega$ is. Meanwhile, $\sqrt{\lambda_{\max}(H_\nu(\omega^*_{opt}))}<\frac{\eta(\omega^*_{aopt})}{\tau}<1$, which implies that SORLopt will converge faster than SORLaopt.

Numerical results for this example are reported in Table~\ref{table2}. From~Table~\ref{table2}, we find that SORLopt always performs better than SORLo and SORLaopt in terms of IT and CPU. SORLopt and SORLno also take the same number of iteration, but SORLopt is superior to SORLno with respect to CPU. In addition, SORLo is better than SORLaopt in terms of IT and CPU. According to IT, CPU and RES, NT is also the best one among the tested methods for this example.

\setlength{\tabcolsep}{5.0pt}
\begin{table}[!h]\small
\centering
\caption{Parameters for Example~\ref{Example4.2}.}\label{table-range}
\begin{tabular}{lllllllll}\hline
$m$ &$\nu$  & $\omega_{o}$ & $\omega_{no}$ &$ \omega^*_{aopt}$ &$\omega^*_{opt}$ & $\sqrt{\lambda_{\max}(H_\nu(\omega^*_{opt}))}$& $\frac{\eta(\omega^*_{aopt})}{\tau}$& $\Omega_1$\\ \hline
  $8$  & $0.2358$ &$1.0671$&$0.991$&$0.8354$& $1$ &$0.3335$ & $0.4309$ & $(0.3994,1.3447)$  \\
  $16$  & $0.2458$ &$1.0704$ &$0.984$&$0.8305$& $1$&$0.3476$ & $0.4438$ & $(0.4003,1.3347)$  \\
  $32$  & $0.2489$ &$1.0714$&$0.992$ &$0.8290$& $1$&$0.3520$ & $0.4478$ & $(0.4005,1.3316)$  \\
  $64$  & $0.2497$ &$1.0717$&$0.983$&$0.8286$ & $1$&$0.3531$ & $0.4488$ & $(0.4006,1.3308)$  \\
   \hline
\end{tabular}
\end{table}

\setlength{\tabcolsep}{2.0pt}% 调整表格列间的宽度
\begin{table}[!h]\small
\centering
\caption{Numerical results for Example~\ref{Example4.2}.}\label{table2}
\begin{tabular}{|c|c|cccc|}\hline
\multirow{2}{*}{Method}  & \multirow{2}{*}{} & \multicolumn{4}{|c|}{$m$}  \\ \cline{3-6}
&       &  $8$ & $16$ & $32$ & $64$ \vextra\\ \hline\hline
\multirow{3}{*}{SORLo}
          &  IT     & $20$ & $21$ & $22$ & $22$\\
          & CPU & $0.0027$ &$0.0049$ &$0.0136$ &$0.0518$\\
          & RES & $3.1565\times 10^{-9}$  & $4.0382\times 10^{-9}$ & $3.3876\times 10^{-9}$ & $7.4133\times 10^{-9}$ \\\hline
\multirow{3}{*}{SORLaopt}
          &  IT     & $23$ & $24$ & $25$ & $26$\\
          & CPU & $0.0034$ &$0.0054$ &$0.0135$ &$0.0534$\\
          & RES & $4.2462\times 10^{-9}$  & $5.3287\times 10^{-9}$ & $4.9212\times 10^{-9}$ & $4.00809\times 10^{-9}$ \\\hline
\multirow{3}{*}{SORLopt}
          &  IT     & $13$ & $14$ & $14$ &$ 15$\\
          & CPU & $0.0022$ &$0.0044$ &$0.0118$ &$0.0507$\\
          & RES & $4.9340\times 10^{-9}$  & $3.0140\times 10^{-9}$ & $6.7716\times 10^{-9}$ & $9.7421\times 10^{-9}$ \\\hline
 \multirow{3}{*}{SORLno}
 &  IT     & $13$ & $14$ & $14$ &$ 15$\\
 & CPU & $1.4552$ &$2.1608$ &$5.2270$ &$20.0681$\\
 & RES & $9.0770\times 10^{-9}$  & $9.4222\times 10^{-9}$ & $9.2093\times 10^{-9}$ & $9.7421\times 10^{-9}$ \\\hline
 \multirow{3}{*}{NT}
 &  IT     & $\textbf{2}$ & $\textbf{2}$ & $\textbf{2}$ &$ \textbf{2}$\\
 & CPU & $\textbf{0.0005}$ &$\textbf{0.0009}$ &$\textbf{0.0024}$ &$\textbf{0.0097}$\\
 & RES & $\boldsymbol{1.7015\times 10^{-14}}$  & $\boldsymbol{3.6599\times 10^{-14}}$ & $\boldsymbol{7.4887\times 10^{-14}}$ & $\boldsymbol{1.6687\times 10^{-13}}$ \\ \hline
%%%%%%%%%%%%%%%%%%%%%%%%%%%%%%%%%%%%%%%%%%%%%%%%%%%%%%%%%%%%%%%%%%%%%%%%%%%%%%%%%%
\end{tabular}
\end{table}

\end{exam}

\begin{exam}[\!\!\cite{davis,kema2017}]\label{Example4.3}
Consider the AVE \eqref{eq:ave} with the matrix $A\in\mathbb{R}^{n\times n}$ arises from six different test problems listed in Table~\ref{table3}. These matrices are sparse and symmetric positive definite and $\|A^{-1}\|<1$. In addition, let $b=Ax^*-|x^*|$ with $x^*=[-1,1,-1,1,\cdots,-1,1]^T\in \mathbb{R}^n$.

Table~\ref{table-range2} displays the parameters for Example~\ref{Example4.3}. It follows from Table~\ref{table-range2} that $\omega^*_{aopt}<\omega^*_{opt}$, $\sqrt{\lambda_{\max}(H_\nu(\omega^*_{opt}))}<\frac{\eta(\omega^*_{aopt})}{\tau}<1$ and the range of $\omega$ becomes smaller as the value of $\nu$ becomes larger. Furthermore, all values of $\omega_{no}$, $\omega^*_{aopt}$ and $\omega^*_{opt}$ belong to $\Omega_1$ or $\Omega_2$ while $\omega_{o}$ does not belong to $\Omega_1$ or $\Omega_2$ for the former three test problems and it belongs to $\Omega_1$ for the later three test problems.

Table~\ref{table4} reports the numerical results for Example~\ref{Example4.3}. From Table~\ref{table4}, for the former three test problems, SORLo is divergent and the reason is that $\omega_{o}$ does not belong to $\Omega_1$ or $\Omega_2$ in these cases. For the test problems $Trefethen\_20b$ and $Trefethen\_200b$, all methods converge within $100$ numbers of iteration. For the $Trefethen\_20000b$ problem, RES of all methods do not reach about $O(10^{-9})$ within $100$ numbers of iteration $($Indeed, it is still the case if we set the maximal iteration number to $2000$ or more$)$.

In conclusion, for the former five test problems in this example, we find that SORLopt is always better than SORLo and SORLaopt in terms of IT and CPU; SORLno is better than SORLopt with respect to IT $($see the possible reason in Remark~\ref{rem:opt}$)$ but behaves worst in terms of CPU; SORLaopt is better than SORLo in terms of IT and CPU and NT is the best one among the tested methods in terms of IT and CPU and RES.

However, for the test problem $Trefethen\_20000b$, if we set $\text{RES} \le 10^{-6}$ as the new stopping criterion, then SORLo, SORLaopt, SORLopt, SORLno and NT converge at iterative steps $53$, $22$, $14$, $12$ and $\textbf{2}$, respectively, and the elapsed CPU time are $12.4365$, $9.4597$, $\textbf{8.6735}$, $13575.4682$ and $9.5410$ seconds, respectively. The RES for SORLo, SORLaopt, SORLopt, SORLno and NT are $7.8509\times 10^{-7}$, $6.7046 \times 10^{-7}$, $4.7848\times 10^{-7}$, $9.6870\times 10^{-7}$ and $3.4955\times 10^{-8}$, respectively. Particularly, in this case, SORLopt is superior to NT in terms of CPU.

%%%%%%%%%%%%%%%%%%%%%%%%%%%%%%%
\setlength{\tabcolsep}{5.0pt}
\begin{table}[!h]\small
\centering
\caption{Test problems for Example~\ref{Example4.3}.}\label{table3}
\begin{tabular}{ll|ll}\hline
 Problem  & $n$ & Problem&  $n$ \\ \hline
  $mesh1e1$  & $48$ & $Trefethen\_20b$ &  $19$ \\
  $mesh1em1$  & $48$ & $Trefethen\_200b$ &  $199$ \\
  $mesh2e1$  & $306$ & $Trefethen\_20000b$ &  $19999$ \\
   \hline
\end{tabular}
\end{table}
%%%%%%%%%%%%%%%%%%%%%%%%%%%%%%%%%%%%%
\setlength{\tabcolsep}{3.0pt}
\begin{table}[!h]\small
\centering
\caption{Parameters for Example~\ref{Example4.3}.}\label{table-range2}
\begin{tabular}{lllllllll}\hline
Problem &$\nu$  &$\omega_{o}$&$\omega_{no}$&$\omega^*_{aopt}$ &$\omega^*_{opt}$ & $\sqrt{\lambda_{\max}(H_\nu(\omega^*_{opt}))}$& $\frac{\eta(\omega^*_{aopt})}{\tau}$ & Range\\ \hline
  $mesh1e1$ & $0.5747$ &$1.2105$&$0.946$&$0.7102$&$0.8218$&$0.7301$ & $0.7588$ & $\Omega_1=(0.4361,1.0753)$  \\
  $mesh1em1$   & $0.6397$&$1.2498$&$0.921$& $0.6929$ & $0.7848$& $0.7845$& $0.8040$ & $\Omega_1=(0.4460,1.0367)$  \\
  $mesh2e1$  & $0.7615$ &$1.3438$&$0.944$&$0.6641$ &$0.7210$& $0.8717$ &  $0.8793$& $\Omega_2=(0.4692,0.9413)$  \\
  $Trefethen\_20b$  &$0.4244$ &$1.1372$&$0.939$&$0.7569$ &$0.9114$ & $0.5783$ &  $0.6365$& $\Omega_1=(0.4175,1.1785)$  \\
  $Trefethen\_200b$ & $0.4265$&$1.1381$ &$0.940$&$0.7561$& $0.9102$& $0.5807$ &  $0.6384$&$\Omega_1=(0.4177,1.1769)$\\
  $Trefethen\_20000b$&  $0.4268$&$1.1382$&$0.943$ &$0.7561$ &$0.9101$ & $0.5810$ &  $0.6387$&$\Omega_1=(0.4177,1.1767)$\\\hline
\end{tabular}
\end{table}

%%%%%%%%%%%%%%%%%%%%%%
\setlength{\tabcolsep}{2pt}
\begin{table}[!h]\small
\centering
\caption{Numerical results for Example \ref{Example4.3} ($\text{RES} \le 10^{-8}$).}\label{table4}
\begin{tabular}{|c|c|ccccc|}\hline
\multirow{2}{*}{Problem}  &  & \multicolumn{5}{|c|}{Method}  \\ \cline{3-7}
                         &  &  $SORLo$   &  $SORLaopt$ & $SORLopt$ & $SORLno$ & $NT$ \vextra\\\cline{3-7}
                         \hline\hline
\multirow{3}{*}{mesh1e1}
                         &  IT          & $-$ & $42$  & $32$ & $23$ & $\textbf{2}$ \\
                         & CPU          & $-$ &$0.0047$   &$0.0039$  &$1.6561$& $\textbf{0.0006}$ \\
                         & RES          & $-$  & $9.5347\times 10^{-9}$ & $9.2561\times 10^{-9}$ & $9.5398\times 10^{-9}$ & $\boldsymbol{9.7580\times 10^{-15}}$\\\hline
\multirow{3}{*}{mesh1em1}
                         &  IT          & $-$ & $43$ & $35$ & $24$ & $\textbf{2}$ \\
                         & CPU          & $-$ &$0.0064$ &$0.0057$ &$1.6590$& $\textbf{0.0006}$ \\
                         & RES          & $-$  & $8.5392\times 10^{-9}$ & $5.9469\times 10^{-9}$ & $9.7687\times 10^{-9}$ & $\boldsymbol{1.7826\times 10^{-14}}$ \\\hline
\multirow{3}{*}{mesh2e1}
                         &  IT          & $-$ & $53$ & $46$ & $26$ & $\textbf{2}$\\
                         & CPU          & $-$ &$0.0252$ &$0.0250$ &$3.0330$& $\textbf{0.0014}$ \\
                         & RES          & $-$  & $7.1097\times 10^{-9}$ & $7.6843\times 10^{-9}$ & $9.9663\times 10^{-9}$ & $\boldsymbol{1.2027\times 10^{-13}}$ \\\hline
\multirow{3}{*}{Trefethen\_20b}
&  IT          & $68$ & $27$ & $18$ & $16$ & $\textbf{2}$ \\
& CPU          & $0.0036$ &$0.0028$ &$0.0013$ &$1.3719$& $\textbf{0.0005}$ \\
& RES          & $7.9973\times 10^{-9}$  & $8.2040 \times 10^{-9}$ & $5.6398 \times 10^{-9}$ & $9.9390\times 10^{-9}$ & $\boldsymbol{4.2074 \times 10^{-14}}$ \\\hline
   \multirow{3}{*}{Trefethen\_200b}
   &  IT          & $69$ & $27$ & $18$ & $16$ & $\textbf{2}$ \\
   & CPU          & $0.0141$ &$0.0114$ &$0.0102$ &$3.8434$& $\textbf{0.0016}$ \\
   & RES          & $9.2407\times 10^{-9}$  & $8.6263\times 10^{-9}$ & $6.1932 \times 10^{-9}$ & $9.4897\times 10^{-9}$ & $\boldsymbol{6.2035\times 10^{-12}}$ \\\hline
 \multirow{3}{*}{Trefethen\_20000b}
   &  IT          & $-$ & $-$ & $-$ & $-$ & $-$ \\
   & CPU          & $-$ &$-$ &$-$ &$-$& $-$ \\
   & RES          & $-$  & $-$ & $-$ & $-$ & $-$\\\hline
\end{tabular}
\end{table}
\end{exam}

As stated in the former three examples, NT performs very well and it is better than all SOR-like methods in most cases, except only for the CPU in the test problem $Trefethen\_20000b$. In addition, SORLo, SORLaopt and SORLno always perform no better than SORLopt in all our tests in terms of CPU. For this reason, in the last example, we will omit reporting numerical results by SORLo, SORLaopt and SORLno. The next example is about the case that $\nu$ is known, thus it is easy to obtain the optimal iteration parameter for SORLopt.

\begin{exam}\label{Example4.4}
The data of the AVE~\eqref{eq:ave} is randomly generated according to~\cite{crfp2016} and the code is free available in {\em http://orizon.mat.ufg.br/admin/pages/11432-codes}. In this example, we run the command $($MATLAB expression and the $data.m$ is download from the above website$)$
$$
``[A,b,\sim,xini,\sim,\sim,normAinv,\sim,\sim,\sim]=data(n,density,cond)"
$$
to generate relatively well-conditioned matrices $A$ with $cond=10$. Since $\nu = normAinv$ is known, the optimal iteration parameter for SORLopt is easy to obtain for this example. ``density" represents the approximate density of matrix $A$. A detailed description of the data can be found in~\cite{crfp2016} or annotations in their codes. In order to avoid the case that both methods are hard to achieve the accuracy of $\text{RES} \le 10^{-8}$ within $100$ numbers of iteration $($like problem $Trefethen\_20000b$ in Example~\ref{Example4.3}$)$, we set $\text{RES} \le 10^{-6}$ as the stopping criterion for this example. In addition, we set $y^{(0)} = x^{(0)}=xini$ for SORLopt. For each set of $(n,density)$, we test ten problems and the averages of IT, CPU and RES are reported in Table~\ref{table7} $($``dens" is short for ``density"$)$. From Table~\ref{table7}, we find that NT always needs less IT than SORLopt. However, SORLopt elapses less CPU than NT, especially for the large sparse AVEs.

\setlength{\tabcolsep}{1.0pt}
\begin{table}[!h]\small
\centering
\caption{Numerical results for Example~\ref{Example4.4}.}\label{table7}
\begin{tabular}{|c|c|ccccc|}\hline
 \multirow{2}{*}{Method}  &  & \multicolumn{5}{|c|}{$n$}  \\ \cline{3-7}
  &     &  $1500 (dens = 1)$ & $1500 (dens = 0.8)$ & $3000 (dens = 0.4)$ & $3000 (dens = 0.1)$ & $10000 (dens = 0.003)$ \vextra\\ \hline\hline
%%%%%%%%%%%%%%%%%%%%%%%%%%%%%%%%%%%%%%%%%%%%%%%%%%%%%%%%%%%%%%%%%%%%%%%%%%%%%%%%%%
\multirow{3}{*}{NT} & IT  & $\textbf{2.9}$ & $\textbf{3.0}$ & $\textbf{3.4}$ & $\textbf{3.0}$& $\textbf{3.3}$ \\
                 & CPU & $0.8996$ & $0.9948$ & $8.7809$ & $2.3885$& $15.0571$ \\
           & RES & $2.7825\times 10^{-9}$ & $9.6605\times 10^{-9}$  & $5.5251\times 10^{-10}$ &  $1.3907\times 10^{-9}$& $2.2440\times 10^{-10}$\\\hline
 \multirow{3}{*}{SORLopt} & IT  & $7.6$ & $6.7$ & $8.6$ & $7.8$& $8.4$\\
                 & CPU & $\textbf{0.6876}$ & $\textbf{0.8100}$ & $\textbf{3.1738}$ & $\textbf{1.0463}$& $\textbf{3.9178}$\\
            & RES &$2.3802\times 10^{-7}$ & $2.3404\times 10^{-7}$ & $2.9646\times 10^{-7}$ & $1.9205\times 10^{-7}$& $1.5401\times 10^{-7}$\\\hline
\end{tabular}
\end{table}

\end{exam}

The numerical results of this section allow us to conclude that SORLopt is better than SORLo and SORLaopt in terms of IT and CPU. SORLopt is superior to SORLno in terms of CPU. However, SORLno is no worse than SORLopt in terms of IT, which demonstrates that our optimal iteration parameter is not exactly the optimal one, it is optimal in the sense that it minimizes $\|T_\nu(\omega)\|_2$, an upper bound of the linear convergence factor of the SOR-like iteration method in the metric $\| |\cdot |\|_{\omega}$ for $(e_{k}^x,e_{k}^y)$. Nevertheless, SORLopt has performed pretty good. In addition, SORLaopt is better than SORLo in some situations. NT is parameter-independent and it is of high efficiency in terms of iteration number. However, in terms of CPU, it will not always be so competitive and sometimes SORLopt is better than it for solving small to middle scale dense or large sparse AVEs. In practice, if $\nu$ is known or can be computed efficiently, SORLopt is a good alternative for solving AVEs. However, in general, even though there exists effective algorithms for computing $\|A^{-1}\|_2$, we cannot necessarily achieve it for large scale problems in an amount of time acceptable to the readers. Thus, the advantages of SORLopt will be weakened or even vanish if the optimal iteration parameter can not be available efficiently. For more detail about the solution of algebraic eigenvalue problems, the reader is referred to \cite{bai2000,wilk1965,cullum1992} and references therein.
%%%%%%%%%%%%%%%%%%%%%%%结论%%%%%%%%%%%%%%%%%%%%%%%%%%%%%
\section{Conclusions}\label{sec:con}
\qquad In this paper, by revisiting the convergence conditions of the SOR-like iteration method proposed in \cite{kema2017} for solving the AVE~\eqref{eq:ave}, the convergent range of the iteration parameter is given and optimal and approximate optimal iteration parameters for the SOR-like iteration method are determined. Our analysis is from the view of the iteration error, which is different from that of iteration matrix \cite{guwl2019}. Furthermore, our optimal and approximate optimal iteration parameters are iteration-independent. Several numerical examples are provided to illustrate that the SOR-like iteration method with our optimal iteration parameter converges faster than that with the approximate optimal parameter proposed in \cite{guwl2019} on solving the AVE \eqref{eq:ave} with $\|A^{-1}\|_2<1$. In addition, our approximate optimal parameter behaves better than that of \cite{guwl2019} in some situations.

Comparing with the generalized Newton method~\cite{mang2009}, the coefficient matrix during the SOR-like iteration is fixed while that of the generalized Newton method is varying. This special intrinsic property allows us to deal with the most expensive implementation, the inverse computation, in a `once-for-all' manner. Hence, the cost per iteration reduces to the matrix-vector productions, which is much lower than the inverse computation. Numerical results show that, in terms of CPU time, the SOR-like iteration method with our optimal iteration parameter sometimes performs better than the generalized Newton method for solving small to middle scale dense or large sparse AVEs.

However, computing the optimal or approximate optimal parameters involved $\|A^{-1}\|_2$, which may be time consuming, especially for large scale problems. In conclusion, if $\|A^{-1}\|_2$ is known or can be efficiently computed, the SOR-like iteration method with our optimal iteration parameter is a good alternative for solving AVE \eqref{eq:ave} with $\|A^{-1}\|_2<1$.

%\section*{Acknowledgements}
%\qquad The first author is grateful to Dr. Yi-Fen Ke and Dr. Chang-Feng Ma from Fujian Normal University for sharing their MATLAB codes.

\section*{\appendix~The proofs of $\boldsymbol{\omega_{i}(\nu) (i=1,\,2,\,3,\,4)}$ being real functions with respect to $\boldsymbol{\nu}$}
\qquad In this section, we prove that $\omega_1(\nu)$ and $\omega_2(\nu)$ are real functions for $\nu \in (0,\frac{\sqrt{2}}{2})$ and $\omega_3(\nu)$ and $\omega_4(\nu)$ are real functions for $\nu \in (\frac{\sqrt{2}}{2}, 1)$, respectively. The proofs are separated in two cases.
\begin{proof}
Let
\begin{align*}
\Delta_1(\nu) &=  -(8\nu^3 - 16\nu^2 + 4\nu -1) + \frac{\nu}{\gamma}(8\nu^3 + 30 \nu^2 - 48\nu + 10),\\
\Delta_2(\nu) &= (8\nu^3 + 16\nu^2 + 4\nu +1) - \frac{\nu}{\zeta}(8\nu^3 - 30 \nu^2 - 48\nu - 10),\\
\Delta_3(\nu) &= -(8\nu^3 - 16\nu^2 + 4\nu -1) - \frac{\nu}{\gamma}(8\nu^3 + 30 \nu^2 - 48\nu + 10),\\
\end{align*}
where
\begin{align*}
\gamma &= \sqrt{-(\nu-1)(\nu+5)},\\
\zeta &= \sqrt{-(\nu+1)(\nu-5)}.
\end{align*}

\begin{itemize}
  \item {\bf Case~1:} when $0 < \nu < \frac{\sqrt{2}}{2}$. Then $\gamma = \sqrt{-(\nu-1)(\nu+5)} > 0 $ and $\zeta = \sqrt{-(\nu+1)(\nu-5)}>0$. Thus, $\omega_1(\nu)$ and $\omega_2(\nu)$ are real functions with respect to $\nu$ if  $\Delta_1(\nu)> 0$ and $\Delta_2(\nu)> 0$, respectively.

      \qquad $\Delta_1(\nu)> 0$ is equivalent to
      \begin{equation}\label{ie:p12}
       p_1(\nu) \doteq 8\nu^4+30\nu^3-48\nu^2+10\nu > (8\nu^3-16\nu^2+4\nu-1)\gamma
       \doteq p_2(\nu)\gamma
      \end{equation}
      Since $p_2'(\nu)=4(6\nu^2-8\nu+1)>0$ in $(0, \frac{4-\sqrt{10}}{6})$ and $p_2'(\nu)<0$ in $(\frac{4-\sqrt{10}}{6}, \frac{\sqrt{2}}{2})$, $p_2$ is strictly monotonously increasing in $(0, \frac{4-\sqrt{10}}{6})$ and strictly monotonously decreasing in $(\frac{4-\sqrt{10}}{6}, \frac{\sqrt{2}}{2})$. Thus $p_2(\nu)\le p_2(\frac{4-\sqrt{10}}{6})=\frac{-83+20\sqrt{10}}{27}<0$. On the other hand, $p_1(\nu)=\nu(4\nu-1)(2\nu+10)(\nu-1)\ge 0$ when $0<\nu \leq \frac{1}{4}$ and $p_1(\nu)<0$ when $\frac{1}{4}<\nu<\frac{\sqrt{2}}{2}$. Thus, from \eqref{ie:p12}, we have $\Delta_1 > 0$ when $0<\nu \leq \frac{1}{4}$. Now we turn to $\frac{1}{4}<\nu<\frac{\sqrt{2}}{2}$, then \eqref{ie:p12} is equivalent to
      \begin{equation*}
       (8\nu^3-16\nu^2+4\nu-1)^2\gamma^2> (8\nu^4+30\nu^3-48\nu^2+10\nu)^2,
      \end{equation*}
      that is
      \begin{align*}
      p_3(\nu)&\doteq (2\nu^2-1)^2(\nu+5)(\nu-1)(-32\nu^2+8\nu-1)\\
       &\doteq (2\nu^2-1)^2(\nu+5)(\nu-1)p_4(\nu) > 0.
      \end{align*}
     Since $p_4(\nu)$ is strictly monotonously decreasing in $(\frac{1}{4},\frac{\sqrt{2}}{2})$, $p_4(\nu)<p_4(\frac{1}{4})= -1<0$. Then we have $p_3(\nu)> 0$. Therefore, if $0<\nu<\frac{\sqrt{2}}{2} $, it holds that $\Delta_1 > 0$.

  \qquad $\Delta_2 > 0$ is equivalent to
  \begin{equation}\label{ie:q12}
  q_1(\nu) \doteq 8\nu^4-30\nu^3-48\nu^2-10\nu < (8\nu^3+16\nu^2+4\nu+1)\zeta
  \doteq q_2(\nu)\zeta.
  \end{equation}
  Since $q_1(\nu)=\nu(\nu-5)(\nu+1)(8\nu+2)$, $q_1(\nu)<0$. On the other hand, $q_2(\nu)>0$. Thus~\eqref{ie:q12} holds.

  \item{\bf Case~2:} when $\frac{\sqrt{2}}{2} < \nu < 1$. Then $\gamma = \sqrt{-(\nu-1)(\nu+5)} > 0 $. Thus, $\omega_3(\nu)$ and $\omega_4(\nu)$ are real functions with respect to $\nu$ if $\Delta_1(\nu)> 0$ and $\Delta_3(\nu) > 0$, respectively.

      \qquad $\Delta_1 > 0$ is equivalent to

  \begin{equation}\label{ie:g12}
  g_1(\nu)\gamma \doteq
  (8\nu^3 - 16\nu^2 + 4\nu -1)\gamma < (8\nu^4 + 30 \nu^3 - 48\nu^2 + 10\nu)\doteq g_2(\nu).
  \end{equation}
  Since $g_1'(\nu)=4(6\nu^2-8\nu+1)<0$ in $(\frac{\sqrt{2}}{2},1)$, $g_1(\nu)$ is strictly monotonously decreasing in $(\frac{\sqrt{2}}{2},1)$. Thus $g_1(\nu)<g_1(\frac{\sqrt{2}}{2})=4\sqrt{2}-9<0$. On the other hand, $ g_2(\nu)=\nu(4\nu-1)(2\nu+10)(\nu-1)$ implies that $g_2(\nu)<0$ whenever $\nu \in (\frac{\sqrt{2}}{2},1)$. Thus, \eqref{ie:g12} is equivalent to
  \begin{equation*}%\label{ie:g13}
  g_3(\nu)\doteq-128\nu^8-480\nu^7+892\nu^6+304\nu^5-776\nu^4+56\nu^3+171\nu^2-44\nu+5 > 0.
  \end{equation*}
  Since
  \begin{align*}
  g_3(\nu)&=(2\nu^2-1)^2(\nu+5)(\nu-1)(-38\nu^2+8\nu-1)\\
  &\doteq(2\nu^2-1)^2(\nu+5)(\nu-1)g_4(\nu),
  \end{align*}
  $2\nu^2-1>0$, $\nu+5>0$, and $\nu-1<0$, it is sufficient to prove $g_4(\nu)< 0$. Indeed, $g_4'(\nu)=-76\nu+8<0$ implies that $g_4(\nu)<g_4(\frac{\sqrt{2}}{2})=-20+4\sqrt{2}<0$.

  \qquad $\Delta_3 > 0$ is equivalent to
  \begin{equation}\label{ie:g14}
  g_5(\nu)\gamma \doteq
  (8\nu^3 - 16\nu^2 + 4\nu -1) \gamma< -(8\nu^4 + 30 \nu^3 - 48\nu^2 + 10\nu)\doteq g_6(\nu).
  \end{equation}
  Earlier, we have proved that $g_5(\nu)=g_1(\nu)<0$. On the other hand, $g_6(\nu)=-\nu(4\nu-1)(2\nu+10)(\nu-1)>0$. Thus, \eqref{ie:g14} holds.
\end{itemize}
\end{proof}

%%%%%%%%%%%%%%%%%%%%%%%%%%参考文献%%%%%%%%%%%%%%%%%%%%%%

\clearpage

\end{document}